\documentclass[12pt, leqno]{article}
\usepackage{amsthm,amsmath,amssymb,bbm,geometry,epsfig,hyperref,enumerate,comment,nicefrac,listings,graphicx,color,relsize}
\numberwithin{equation}{section}
\newtheorem{lemma}{Lemma}[section]

\newtheorem{algo}[lemma]{Framework}

\providecommand{\N}{{\ensuremath{\mathbb{N}}}}

\providecommand{\R}{{\ensuremath{\mathbb{R}}}}
\providecommand{\B}{\mathcal{B}}

\renewcommand{\P}{\mathbbm{P}}

\providecommand{\bS}{\mathbb{S}}

\renewcommand{\S}{\mathcal{S}}
\providecommand{\E}{{\ensuremath{\mathbbm{E}}}}

\providecommand{\N}{{\ensuremath{\mathbbm{N}}}}

\providecommand{\R}{{\ensuremath{\mathbbm{R}}}}

\providecommand{\E}{{\ensuremath{\mathbb{E}}}}
\newcommand{\F}{{\ensuremath{\mathcal{F}}}}

\newcommand{\Y}{{\ensuremath{\mathcal{Y}}}}

\DeclareFontEncoding{LS1}{}{}
\DeclareFontSubstitution{LS1}{stix}{m}{n}
\DeclareMathAlphabet{\mathscr}{LS1}{stixscr}{m}{n}
\begin{document}

\title{%\vspace{-3.9cm}
	 Deep splitting method for  parabolic PDEs
	%Deep splitting method: Numerical approximations\\ for high-dimensional partial differential equations\\ using deep learning and the splitting-up method
}

\author{Christian Beck$^1$, Sebastian Becker$^2$, 
	Patrick Cheridito$^3$, \\
	 Arnulf Jentzen$^4$, and Ariel Neufeld$^5$
	\bigskip
	\\
	\small{$^1$ Department of Mathematics, ETH Zurich,}\\
	\small{Switzerland, e-mail:  
		christian.beck@math.ethz.ch}
	\smallskip
	\\
	\small{$^2$ %RiskLab, 
		Department of Mathematics, ETH Zurich,}\\
	\small{Switzerland, e-mail:  
		sebastian.becker@math.ethz.ch}
	\smallskip
	\\
	\small{$^3$ %RiskLab, 
		Department of Mathematics, ETH Zurich, }\\
	\small{Switzerland, e-mail:   patrick.cheridito@math.ethz.ch}
	\smallskip
	\\
	\small{$^4$ Institute for Analysis and Numerics, Faculty of Mathematics and Computer Science, }\\
	\small{University of Münster, e-mail:  ajentzen@uni-muenster.de}
	\smallskip
	\\
	\small{$^5$ Division of Mathematical Sciences, School of Physical and Mathematical Sciences,}
	\\
	\small{Nanyang Technological University, Singapore,
		e-mail: ariel.neufeld@ntu.edu.sg}
}

\date{April 2021}
%%%%%%%%%%%%%%%%%%%

%%%%%%%%%%%%%%%5
\maketitle
\vspace{-0.6cm}
\begin{abstract}
In this paper we introduce a numerical method for nonlinear parabolic PDEs that combines operator 
splitting with deep learning. It divides the PDE approximation problem into a sequence of separate learning 
problems. Since the computational graph for each of the subproblems is comparatively small, the approach 
can handle extremely high-dimensional PDEs. We test the method on different examples from physics, 
stochastic control and mathematical finance. In all cases, it yields very good results in up to 10,000 
dimensions with short run times.\\[2mm]
\end{abstract}

\newpage

%\tableofcontents
%
%

%
%
\section{Introduction}
In this paper we derive a numerical scheme for parabolic partial differential equations (PDEs) of the form
\begin{equation} \label{eq:pde}
\begin{aligned}
\tfrac{\partial}{\partial t}u(t,x) &=  F \big(x, u(t,x),  \nabla_x u ( t,x ) \big) + \tfrac{1}{2}
\operatorname{Trace}\!\big( 
\sigma(x) \sigma^*(x) 
\operatorname{Hess}_x u( t,x )\big), 
%	\\u(0,x) &= \varphi(x),
\end{aligned}
\end{equation}
$(t,x) \in (0,T] \times \R^d$, with initial condition $u(0,x) = \varphi(x)$,
where $F \colon \mathbb{R}^d \times \mathbb{R} \times \mathbb{R}^d \to \mathbb{R}$, 
$\sigma \colon \mathbb{R}^d \to \mathbb{R}^{d \times d}$, and $\varphi \colon \mathbb{R}^d \to \mathbb{R}$ 
are appropriate continuous functions. %and $\operatorname{Trace}$ denotes the trace of a $d$-dimensional square matrix. 
Such PDEs describe various phenomena in nature, engineering, 
economics, and finance. They typically do not admit closed form solutions and, therefore, have to be solved numerically. 
In some applications, the dimension $d$ can be high. For instance, in physics and engineering 
applications, $x \in \mathbb{R}^d$ typically models the coordinates of all components of a given system,
whereas in derivative pricing and optimal investment problems, $d$ usually corresponds to the number of underlying assets. 
Many classical PDEs, such as the standard heat and Black--Scholes equations
are linear. Using the Feynman--Kac representation, their solutions can efficiently 
be approximated in high dimensions with simple Monte Carlo averages.
But if constraints or frictions are taken into account, or the PDE describes 
a control problem, the function $F$ is no longer linear 
and equation \eqref{eq:pde} becomes much more challenging to solve for large $d$.

Numerical methods for PDEs have a long history. Classical approaches like finite differences and finite elements (see, e.g., \cite{Braess2007FEM,Thomee2003PDEsAndNumerics,Thomee1997}) are deterministic.
In their standard form, they work well for $d=1,2$ and $3$, but their complexity grows exponentially in $d$. 
To tackle higher dimensional problems, different simulation-based approaches have been developed
that exploit a stochastic representation of the solution of the PDE. For instance, \cite{BallyPages2003,
	BenderDenk2007, 
	Bender2015Primal, 
	BouchardTouzi2004,
	Chassagneux2014, 
	ChassagneuxCrisan2014,
	ChassagneuxRichou2015, 
	ChassagneuxRichou2016,
	CrisanManolarakis2010,
	CrisanManolarakis2012, 
	CrisanManolarakis2014, 
	CrisanManolarakisTouzi2010, 
	DelarueMenozzi2006,
	DouglasMaProtter, 
	GobetLabart2010, 
	GobetLemor2008Numerical, 
	GobetLemorWarin2005,  
	GobetLopezSalasTurkedjiev2016,
	GobetTurkedjiev2016,
	GobetTurkedjiev2016MathComp,  
	HuijskensRuijterOosterlee2016,
	LabartLelong2013,
	LemorGobetWarin2006, 
	LionnetDosReisSzpruch2015, 	      
	MaProtterSanMartin2002, 
	MaProtterYong1994, 
	MaYong1999, 
	MilsteinTretyakov2006, 
	MilsteinTretyakov2007, 
	Pham2015,
	RuijterOosterlee2015,RuijterOosterlee2016,Ruszczynski2017Dual,Turkedjiev2015,Zhang2004}
use BSDE representations of PDEs and study approximation methods based on recursive polynomial regressions, 
\cite{Labordere2012, 
	Labordereetal2016arxiv,
	LabordereTanTouzi2014,
	McKean1975,
	SkorohodBranchingDiffusion1964, 
	Watanabe1965Branching} investigate methods based on branching diffusion processes, 
and ~\cite{LinearScaling, MultilevelPicard, hutzenthaler2018overcoming,hutzenthaler2019overcoming,HutzenthalerKruse17} analyze full-history recursive multilevel Picard methods.
%Recently, a series of numerical approximation methods for high-dimensional PDEs have been 
%developed based on the idea to reformulate the PDE as a stochastic learning problem.
%This makes them accessible to the toolkit of deep learning; see [E, Han, Jentzen etc.].
%Except of MLP approximation methods, most of the above named approximation methods are only 
%applicable in the case where the PDE/BSDE dimension $ d $ is rather 
%small or work exclusively in the case of serious restrictions on the 
%parameters or the type of the considered PDE.
%small terminal/initial conditions). 
%
%
Recently, numerical methods for high-dimensional PDEs based on the idea to reformulate the PDE as a stochastic learning problem have been proposed in \cite{EHanJentzen2017,HanEJentzen17}. This
opens the door to the application of deep learning; see, e.g.,
\cite{DeepKolmogorov,BeckEJentzen17,becker2018deep,
berg2018unified,chan2018machine, EYu17,FarahmandNabiNikovski17, FujiiTakahashiTakahashi17,goudenege2019machine, han2018convergence,HenryLabordere17,hure2019some,jacquier2019deep,LongLuMaDong17,lye2019deep,magill2018neural,Raissi18,SirignanoDGM2017} for modifications and extensions. There are also already a few papers studying the convergence of deep learning based approximation methods for PDEs. For instance, \cite{han2018convergence,SirignanoDGM2017} derive convergence results without information on the convergence speed, whereas \cite{berner2018analysis,elbrachter2018dnn,grohs2018proof,hutzenthaler2019proof,jentzen2018proof,kutyniok2019theoretical} provide convergence and tractability results with dimension-independent convergence rates and error constants depending polynomially on the dimension.
%%%%%%%%%%%%%%%%%%%%5 

In this paper we develop a new deep learning method for parabolic PDEs that splits the 
differential operator into a linear and a nonlinear part. More precisely, we write
\begin{equation}
F \big(x, u(t,x), \nabla_x u ( t,x ) \big) 
= 
\langle \mu(x), \nabla_x u ( t,x ) \rangle_{\R^d} 
+ f\big(x, u(t,x), \nabla_x u ( t,x )\big)
\end{equation}
for suitable continuous functions $\mu \colon \mathbb{R}^d \to \mathbb{R}^d$ and 
$f \colon \mathbb{R}^d \times \mathbb{R} \times \mathbb{R}^d \to \mathbb{R}$. This decomposition is not unique.
But the idea is that $\mu$ is chosen such that the nonlinearity $f \big(x, u(t,x), \nabla_x u ( t,x ) \big)$ becomes small. 
Then we solve the PDE iteratively over small time intervals by approximating $f \big(x, u(t,x),  \nabla_x u ( t,x ) \big)$
and using the Feynman--Kac representation locally. This requires a recursive computation of 
conditional expectations. We approximate them by formulating them as minimization problems 
that can be approached with deep learning. This decomposes the PDE approximation problem into 
a sequence of separate learning problems. Since the computational graph for each of the 
subproblems is comparatively small, the method works for very high-dimensional problems.

%%%%%%%%%%%%%%%%%%%%%%%%%%%%%%%%%%%%%%%%%%%%%%
The rest of the paper is organized as follows. In Section \ref{sec:derivation} we introduce the 
framework and derive the deep splitting method. In Section \ref{sec:examples} we test the approach 
on five different classes of high-dimensional nonlinear 
PDEs: Hamilton--Jacobi--Bellman (HJB) equations, nonlinear Black--Scholes equations, 
Allen--Cahn-type equations, nonlinear heat equations, 
and sine-Gordon-type equations.
%%%%%%%%%%%%%%%%%%%%%%%%%%%%%%%%%%%%%%%%%%%%%%%%%%

%%%%%%%%%%%%%%%%%%%%%%%%%%%%%%%%%%%%%%%%%%%%%%%%%%%%%%%%%%%%%5

%\section{Derivation for a stochastic heat equation}
%\input{Special_Case.tex}
%%%%%%%%%%%%%%%%%%%%%%%%%%%%%%%%%%%%%%%%%%%%%%%%%%%%%%%%%%%%%%%%%%%%%%%%%%%
\section{Derivation of the proposed approximation algorithm}
\label{sec:derivation}
%\input{General_Case.tex}
%%%%%%%%%%%%%%%%%%%%%%%%%%%%
Fix $ T \in (0,\infty)$ and $d \in \N $. Consider two at most polynomially growing continuous functions
$\varphi \colon \R^d \to \R$ and $f \colon \R^d \times \R \times \R^d \to \R$ together with two Lipschitz continuous functions 
$ \mu \colon \R^d \to \R^d$ and $\sigma \colon  \R^d \to \R^{ d \times d }$.
Assume $u \colon [0,T] \times \R^d \to \R$ is an at most polynomially growing continuous function 
that is $C^{1,2}$ on $(0,T] \times \R^d$ and satisfies the PDE 
\begin{equation}
\begin{split}
\label{eq:defPDE}
\tfrac{\partial}{\partial t}u(t,x)
&=  f\big(x, u(t,x), \nabla_x u (t,x ) \big)
+
\big\langle \mu(x), \nabla_x u ( t,x ) \big\rangle_{ \R^d }\\
& \quad + \tfrac{ 1 }{ 2 }
\operatorname{Trace}\!\big( 
\sigma(x) \sigma^*(x)
\operatorname{Hess}_x u( t,x )
\big), 
\end{split}
\end{equation}
$(t,x) \in (0,T] \times \R^d$, with initial condition $u(0,x)=\varphi(x)$, $x \in \R^d$.
\subsection{Temporal discretization}
\label{subsec:temp-discret}

To approximate the solution $u$ of the PDE \eqref{eq:defPDE}, we discretize the equation 
in time and use a splitting-up method 
(see, e.g., \cite{
	GrekschLisei_ApproximationOfStochasticNonlinearEquationsBySplittingMethod2013,
	GyoengyKrylov_OnTheRateOfConvergenceOfSplittingUpApproximationsForSPDEs2003,
	GyoengyKrylov_OnTheSplittingUpMethodForSPDEs}) 
to obtain a semi-discrete approximation problem. To do this, we choose
$N \in \N$ and let $t_0, t_1, \ldots, t_N\in [0,T]$ be real numbers such that 
\begin{equation}
\label{eq:time-step-discrete}
0 = t_0  < t_1 < \ldots < t_N = T.
\end{equation}
Under appropriate integrability assumptions, it follows from \eqref{eq:defPDE} that for every $t\in[0,T]$, $x\in \R^d$
we have 
\begin{equation}\label{eq:defIntPDE}
\begin{split}
u(t,x)
&  =
\varphi( x ) 
+
\int_{ 0 }^{ t }
f\big( 
x, u(s,x), \nabla_x u ( s,x ) 
\big)
\, ds
%  +
%  \int_{ 0 }^{ t }
%  b\big( x, X_s( x ), ( \nabla X_s )( x ) \big)
%  \, dZ_s(x)
\\
& \quad
+
\int_{ 0 }^{ t }
\Big[
\tfrac{ 1 }{ 2 }
\operatorname{Trace}\!\big( 
\sigma(x) \sigma^*(x) 
\operatorname{Hess}_x u( s,x )
\big)
+
\big\langle \mu(x), \nabla_x u ( s,x ) \big\rangle_{ \R^d }
\Big]
\, ds.
\end{split}
\end{equation}
In particular, for all $n\in\{0,1,\ldots,N-1\}$, $t\in [t_n,t_{n+1}]$ and $x\in\R^d$,
\begin{equation}
\begin{split}
u(t,x) 
& 
= u(t_n,x) 
+
\int_{ t_n }^{ t }
f\big( 
x, u(s,x), ( \nabla_x u )( s,x ) 
\big)
\, ds
%  +
%  \int_{ t_n }^{ t }
%  b\big( x, X_s( x ), ( \nabla X_s )( x ) \big)
%  \, dZ_s(x)
\\
&
\quad  + \int_{t_n}^t 
\Big[
\tfrac{ 1 }{ 2 }
\operatorname{Trace}\!\big( 
\sigma(x) \sigma^*(x) 
\operatorname{Hess}_x u( s, x )
\big)
+
\big\langle \mu(x), \nabla_x u ( s,x ) \big\rangle_{ \R^d }
\Big]
\, ds,
\end{split}
\end{equation}
and therefore, 
\begin{equation}\label{eq:Xn_discretization}
\begin{split}
u(t,x) %X_t (x) 
& \approx 
u(t_n,x) %X_{t_n}(x)
+
\int_{ t_n }^{ t_{ n + 1 } }
f\big( 
x, u(t_n,x), \nabla_x u ( t_n,x ) 
\big)
\, ds
%    +
%    \int_{ t_n }^{ t_{ n + 1 } }
%    b\big( x, X_{t_n}( x ), ( \nabla X_{t_{n}} )( x ) \big)
%    \, dZ_s(x)
\\
& \quad + \int_{t_{n}}^t 
\Big[
\tfrac{ 1 }{ 2 }
\operatorname{Trace}\!\big( 
\sigma(x) \sigma^*(x) 
\operatorname{Hess}_x u( s,x )
\big) 
+ 
\big\langle \mu(x), \nabla_x u( s,x ) \big\rangle_{ \R^d }\Big]\,ds,
\end{split}
\end{equation}
which can be written as 
\begin{equation}\label{eq:Xapproximation}
\begin{split} 
u(t,x)   
&
\approx
u(t_n,x) 
+
f
\bigl(
x,u(t_n,x), \nabla_x u(t_n,x)
\bigr)\,(t_{n+1}-t_n)
\\[1ex]
%& 
%  + 
%   b\big(x, X_{t_n}(x), (\nabla X_{t_n})(x)\big)\,\big(Z_{t_{n+1}}(x)-Z_{t_n}(x)\big)
%  \\
& \quad
+ 
\int_{t_{n}}^t 
\Big[
\tfrac{ 1 }{ 2 }
\operatorname{Trace}\!\big( 
\sigma(x) \sigma^*(x)
\operatorname{Hess}_x u( s,x )
\big) 
+ 
\big\langle \mu(x), \nabla_x u( s,x ) \big\rangle_{ \R^d }\Big]\,ds .
\end{split}
\end{equation}
To derive the splitting-up approximation, we make a few simplifying assumptions, not all of which are 
needed to implement the resulting algorithm. First, we suppose that $\varphi$ has an at most polynomially 
growing gradient $\nabla \varphi \colon \R^d \to \R^d$. Moreover, we assume that there exists a function
$v \colon [0,T] \times\R^d \to \R$ satisfying $v(0, x) = \varphi(x)$, $x \in \R^d$, such that for every $n \in \{0,1,\ldots,N-1\}$,
$v|_{(t_n,t_{n+1}]\times\R^d}$ belongs to $C^{1,2}((t_n,t_{n+1}]\times\R^d,\R)$ with at most polynomially growing 
partial derivatives and for all $n \in \{0,1,\ldots,N-1\}$, $t \in (t_n, t_{n+1}]$ and $x \in \R^d$,
%the function $v$ satisfies
%
\begin{equation}\label{eq:mildFormulationVPDE}
\begin{split}
v(t,x)
& =
v(t_n,x) 
+ 
f \big(x,v(t_n,x),\nabla_x v(t_n,x)\big)\,(t_{n+1}-t_n) \\[1ex]
% & + 
%b(x,V_{t_n}(x),(\nabla V_{t_n})(x))\,(Z_{t_{n+1}}(x)-Z_{t_n}(x)) \\[1ex]
&\quad  + 
\int_{ t_n }^{ t }
\Big[
\tfrac{ 1 }{ 2 }
\operatorname{Trace}\!\big( 
\sigma(x) \sigma^*(x)
\operatorname{Hess}_x v(s,x )
\big)
+
\big\langle \mu(x), \nabla_x v(s, x ) \big\rangle_{ \R^d }
\Big]
\, ds;
\end{split}
\end{equation}
see, e.g.,~Hairer et al.\ \cite[Section 4.4]{HairerHutzenthalerJentzen_LossOfRegularity2015}, 
Deck \& Kruse \cite{DeckKruse_ParametrixMethod2002},
Krylov \cite[Chapter 8]{Krylov_LecturesHoelder1996},
and Krylov \cite[Theorem 4.32]{Krylov_KolmogorovsEquations1998}
for existence, uniqueness, and regularity results for equations of the form
\eqref{eq:mildFormulationVPDE}.
Comparing \eqref{eq:mildFormulationVPDE} to \eqref{eq:Xapproximation} suggests that 
\begin{equation} \label{eq:VapproxX}
v(t_n,x) \approx u(t_n,x) \quad \mbox{for all } n \in \{ 1, \dots, N\}.
\end{equation}
So $v$ is a specific splitting-up type approximation of the function $u$; see, e.g.,~\cite{Dorsek12,
	GrekschLisei_ApproximationOfStochasticNonlinearEquationsBySplittingMethod2013,
	GyoengyKrylov_OnTheRateOfConvergenceOfSplittingUpApproximationsForSPDEs2003,
	GyoengyKrylov_OnTheSplittingUpMethodForSPDEs}.
%
%
%%%%%%%%%%%%%%%%%%
\subsection{An approximate Feynman--Kac representation}

In the next step we derive a Feynman--Kac representation of $v$; see, e.g., Milstein \&  Tretyakov \cite[Section 2]{MilsteinTretyakov_SolvingPSPDEsViaAveraging2009}. Let $B \colon [0,T]\times\Omega\to\R^d$ be a
standard $(\mathcal{F}_t)_{t\in [0,T]}$-Brownian motion on a filtered probability space
$(\Omega, \mathcal{F}, ( \mathcal{F}_t )_{ t \in [0,T] }, \P)$ satisfying the usual conditions.
Consider an $\mathcal{F}_0$/$\B(\R^d)$-measurable random variable
$\xi \colon \Omega\to\R^d$ satisfying $\E[\|\xi\|_{\R^d}^p] < \infty$ for every $p\in (0,\infty)$,
and let $Y\colon [0, T]\times\Omega\to\R^d$ be an $(\mathcal{F}_t)_{t\in [0,T]}$-adapted 
process with continuous sample paths satisfying for every $t\in [0,T]$,
\begin{equation}\label{eq:SDEForY}
Y_{t} = \xi + \int_0^t \mu(Y_s)\,ds + \int_0^t \sigma(Y_s)\,dB_s \quad \mbox{$\P$-a.s.}
\end{equation}
The assumption that $\E[\|\xi\|_{\R^d}^p] < \infty$ for all $p\in (0,\infty)$ 
ensures that 
\begin{equation}
\label{eq:Y-integ-cond}
\sup_{t \in [0,T]} \E\big[\|Y_t\|_{\R^d}^p\big] < \infty
\quad \mbox{for all } p \in (0, \infty);
\end{equation}
see, e.g., Stroock~\cite[Section 1.2]{Stroock1982TopicsInSDEs}.
Moreover, \eqref{eq:mildFormulationVPDE} implies that for all
$n\in\{0,1,\ldots,N-1\}$, $t\in (t_{n},t_{n+1})$ and $x\in\R^d$, one has 
\begin{equation} \label{eq:vt}
\tfrac{\partial}{\partial t} v(t,x) 
= \big\langle \mu(x), \nabla_x v(t, x ) \big\rangle_{ \R^d } 
+ \tfrac{ 1 }{ 2 }
\operatorname{Trace}\!\big( 
\sigma(x) \sigma^*(x) 
\operatorname{Hess}_x v(t,x )
\big)
\end{equation}
from which it follows that for all 
$n\in\{0,1,\ldots,N-1\}$, 
$t\in (T-t_{n+1},T-t_{n})$ and
$x\in\R^d$,
\begin{multline} \label{eq:backwardsEquationPoppingUpInIto}
% \begin{split}
\tfrac{\partial}{\partial t} v(T-t,x) 
+ \big\langle \mu(x), \nabla_x v(T-t, x ) \big\rangle_{ \R^d } 
+ \tfrac{ 1 }{ 2 }
\operatorname{Trace}\big( 
\sigma(x ) \sigma^*(x )
\operatorname{Hess}_x v(T-t, x )
\big)
= 0. 
% \end{split}
\end{multline}
Since for every $n\in\{0,1,\ldots,N-1\}$, $v|_{(t_{n},t_{n+1}]\times\R^d}$ is in $C^{1,2}((t_{n},t_{n+1}]\times\R^d,\R)$, 
we obtain from It\^o's formula that for all $n\in\{0,1,\ldots,N-1\}$ and $t \in [T-t_{n+1},T-t_{n})$ we have 
\begin{equation}
\begin{split}
& v(T-t,Y_t) = v(T- t_{n+1},Y_{t_{n+1}}) + 
\int_{T - t_{n+1}}^{t} 
\big\langle \nabla_x v(T-s,Y_s), 
\sigma(Y_s)\,dB_s  \big\rangle_{\R^d}\\
& + 
\int_{T - t_{n+1}}^{t} 
\tfrac{\partial}{\partial s} v(T-s, Y_s)\,ds
 +
\int_{T - t_{n+1}}^{t} 
\tfrac{ 1 }{ 2 }
\operatorname{Trace}\big( 
\sigma( Y_s ) \sigma^*(Y_s ) 
\operatorname{Hess}_x v(T-s, Y_s )
\big)\,ds \\
& + 
\int_{T - t_{n+1}}^{t}
\big\langle \mu(Y_s ), 
\nabla_x v(T-s, Y_s ) 
\big\rangle_{ \R^d }
\, ds \quad \mbox{$\P$-a.s.},
\end{split}
\end{equation}
which by \eqref{eq:backwardsEquationPoppingUpInIto}, gives that $\P$-a.s.,
\begin{equation}\label{eq:ItoOnOpenIntervalAfterCancellation}
v(T-t,Y_t) 
= 
v(T - t_{n+1},Y_{t_{n+1}})
+ 
\int_{T - t_{n+1}}^t 
\big\langle \nabla_x v(T-s,Y_s), 
\sigma(Y_s)\,dB_s  \big\rangle_{\R^d}. 
\end{equation}
Moreover, since, by assumption, $\sigma \colon \R^d \to \R^{d\times d}$ is Lipschitz continuous and 
$\nabla_x v(t,x)$ at most polynomially growing in $(t,x) \in (t_n, t_{n+1}] \times \R^d$, one obtains from 
\eqref{eq:Y-integ-cond} that for all $n \in \{0,1,\dots,N-1\}$ and $t\in [T-t_{n+1},T-t_{n})$,
\begin{equation}
\int_{T-t_{n+1}}^t \E \Big[\big\| \sigma^*(Y_s) \nabla_x v(T-s,Y_s)\big\|_{\R^d}^2\Big]\, ds <\infty,
\end{equation}  
from which it follows that 
\begin{equation}\label{eq:conditionalExpectation1stTerm}
\E\bigg[\int_{T-t_{n+1}}^{t} 
\big\langle \nabla_x v(T-s,Y_s), 
\sigma(Y_s)\,dB_{s} \big\rangle_{\R^d}\!~\Big|\!~\mathcal{F}_{T-t_{n+1}} \bigg] = 0  \quad \mbox{$\P$-a.s.}
\end{equation}
Together with \eqref{eq:ItoOnOpenIntervalAfterCancellation}, this shows that for all
$n\in\{0,1,\ldots,N-1\}$ and $t \in [T-t_{n+1}, T-t_{n})$ we have
\begin{equation}
\E\Big[
v(T-t,Y_t)\!~\big|\!~\mathcal{F}_{T-t_{n+1}}
\Big] =
\E\Big[
v(t_{n+1},Y_{T-t_{n+1}})\!~\big|\!~\mathcal{F}_{T-t_{n+1}} \Big] \quad \mbox{$\P$-a.s.}
\end{equation}
Since $Y_{T- t_{n+1}}$ is $\mathcal{F}_{T-t_{n+1}}$/$\B(\R)$-measurable, one obtains from the 
tower property of conditional expectations that for all $n\in\{0,1,\ldots,N-1\}$ and $t \in [T-t_{n+1}, T-t_{n})$,
\begin{equation} \label{eq:vcond}
\E\Big[ v(T-t,Y_t) \!~\big|\!~Y_{T-t_{n+1}}\Big] = \E\Big[v(t_{n+1},Y_{T-t_{n+1}}) \!~\big|\!~Y_{T-t_{n+1}}\Big] = v(t_{n+1},Y_{T-t_{n+1}})
\quad \mbox{$\P$-a.s.}
\end{equation}
Since $(Y_t)_{t \in [0,T]}$ has continuous sample paths and
$(v(t,x))_{(t,x)\in (t_n,t_{n+1}]\times\R^d}$ at most polynomially growing first order 
partial derivatives, it follows from \eqref{eq:mildFormulationVPDE} and \eqref{eq:vt} that 
we have for all $n \in \{0,1,\dots, N-1\}$ and $\omega \in \Omega$,
\begin{equation}\label{eq:ptwConv}
\begin{split}
% \begin{multline}
&\lim_{t\nearrow T-t_n} v(T-t,Y_t(\omega)) =
v(t_n,Y_{T-t_n}(\omega)) \\& \qquad +
f\big(Y_{T-t_n}(\omega),v(t_n,Y_{T-t_n}(\omega)), \nabla_x v(t_n,Y_{T-t_n}(\omega))\big)\,(t_{n+1}-t_n).
% \end{multline}
\end{split}
\end{equation}
In addition, since $\sup_{t\in [0,T]} \E[\|Y_t\|_{\R^d}^p]$ $<\infty$ for every $p \in (0, \infty)$, we obtain 
\begin{equation}
\sup_{t \in (t_n, t_{n+1}] }
\E\big[ 
|v(T-t,Y_t)|^{p}
\big] < \infty. 
\end{equation}
So it follows from \eqref{eq:vcond} and \eqref{eq:ptwConv} that for all $n\in\{0,1,\ldots,N-1\}$,
\begin{equation}\label{eq:conditionalExpRZ}
\begin{split}
&v(t_{n+1},Y_{T-t_{n+1}}) = \E\Big[\lim_{t\nearrow T-t_n} v(T-t,Y_t) \!~\big|\!~Y_{T-t_{n+1}}\Big]\\
& =
\E\Bigl[ 
v(t_n,Y_{T-t_n}) +
f\big(Y_{T-t_n}, v(t_n,Y_{T-t_n}), \nabla_x v(t_n,Y_{T-t_n})\big)
\,(t_{n+1}-t_n) 
\!~\big|\!~Y_{T-t_{n+1}} 
\Bigr] \quad \mbox{$\P$-a.s.,}
\end{split}
\end{equation}
which is the Feynman--Kac type representation we were aiming for. Note that the nonlinearity $f$ 
as well as the coefficient functions $\mu$ and $\sigma$ in \eqref{eq:defPDE} do not depend on $t$.
But the above derivation could be extended to time-dependent $f$, $\mu$, and $\sigma$
since the Feynman--Kac formula still holds in this case.

\subsection{Formulation as recursive minimization problems}

We now reformulate \eqref{eq:conditionalExpRZ} as recursive minimization problems. 
It follows from our assumptions that for every $n\in\{1,2,\ldots,N\}$, 
$v(t_{n-1},x) + f(x,v(t_{n-1},x),\nabla_x v(t_{n-1},x))(t_{n}-t_{n-1})$
is at most polynomially growing in $x \in \R^d$. Therefore, we obtain from \eqref{eq:Y-integ-cond} that
\begin{equation}\label{eq:VtnIsInL2}
\E\Big[\big|v(t_{n-1},Y_{T-t_{n-1}}) + f\big(Y_{T-t_{n-1}},
v(t_{n-1},Y_{T - t_{n-1}}), \nabla_x v(t_{n-1},Y_{T-t_{n-1}}\big)\,(t_{n}-t_{n-1})
\big|^2\Big] < \infty.
\end{equation}
%  \end{split}
% \end{equation}
Since $v(t_n, x)$ is continuous in $x \in \R^d$, it follows from the 
factorization lemma and the $L^2$-minimality property of conditional 
expectations (see, e.g., Klenke~\cite[Corollary 8.17]{Klenke_2014})
that for every $n\in\{1,2,\ldots,N\}$ we have 
\begin{multline}\label{eq:formulationAsMinimizationProblem-OLD}
%\begin{split} 
(v(t_{n},x))_{x\in\operatorname{supp}(Y_{T-t_{n}}(\P))}
=  
\operatornamewithlimits{argmin}_{w \in C(\operatorname{supp}(Y_{T-t_{n}}(\P)),\R)}
\E
\Big[ 
\big|w(Y_{T-t_{n}}) 
- \big[
v(t_{n-1},Y_{T-t_{n-1}})\\ 
%
% \quad %\quad\quad\quad\quad
+ f \big(Y_{T-t_{n-1}}, v(t_{n-1},Y_{T-t_{n-1}}), \nabla_x v(t_{n-1},Y_{T-t_{n-1}})\big)\,(t_{n}-t_{n-1})
\big]
\big|^2
\Big]. 
% \end{split}
\end{multline}

\subsection{Deep artificial neural network approximations}
To tackle the minimization problems \eqref{eq:formulationAsMinimizationProblem-OLD} numerically, we 
approximate the functions $v(t_n, .)$, $n \in \{1,2, \dots,N\}$, with neural networks $V_n$. 
More precisely, we choose $k \in \{3,4, ...\}$ and $l \in \N$. Then, we set $\nu = (1+ kl - l)(l+1) + l(d+1)$
and consider functions $V_n \colon \R^{\nu} \times \R^d \to \R$, $n\in\{0,1,\dots,N\}$, 
such that for every $(\theta, x) \in \R^{\nu} \times \R^d$, $V_0$ is given by 
\begin{equation} \label{eq:V0}
V_0(\theta, x) = \varphi(x)
\end{equation}
and for $n \in \{1,2, \dots, N\}$, $V_n$ is of the form 
\begin{equation} \label{eq:Vn}
V_{ n }(\theta,x) = A^{ \theta, (k-1)l(l+1) + l(d+1) }_{l, 1 } 
\circ 
\mathcal{L}_l
\circ
A^{ \theta, (k-2)l(l+1) + l(d+1) }_{l,l } 
%   & \quad
\circ
\ldots
%   \circ
%   \mathcal{L}_k
%   \circ
%   A^{ \nu, \theta, 3 k(k+1) }_{ k, k } 
\circ 
\mathcal{L}_l
\circ 
A^{ \theta, l(d+1) }_{l,l } 
\circ 
\mathcal{L}_l
\circ 
A^{ \theta, 0 }_{d,l }(x),
%  \end{split}
\end{equation}
where for $ r \in \N_0 = \{0\} \cup \N $ and $ i,j\in \N $
with 
$
r + j(i + 1 ) \leq \nu,
$
$ A^{ \theta, r }_{ i, j } \colon \R^i \to \R^j $ is the affine function defined by 
\begin{equation}
A^{ \theta, r }_{i,j }( x )
=
\left(
\begin{array}{cccc}
\theta_{ r + 1 }
&
\theta_{ r + 2 }
&
\dots
&
\theta_{ r + i }
\\
\theta_{ r + i + 1 }
&
\theta_{ r + i + 2 }
&
\dots
&
\theta_{ r + 2 i }
\\
\theta_{ r + 2 i + 1 }
&
\theta_{ r + 2 i + 2 }
&
\dots
&
\theta_{ r + 3 i }
\\
\vdots
&
\vdots
&
\vdots
&
\vdots
\\
\theta_{ r + i( j - 1 ) + 1 }
&
\theta_{ r + i( j - 1 ) + 2 }
&
\dots
&
\theta_{ r +  ij }
\end{array}
\right)
\left(
\begin{array}{c}
x_1
\\
x_2
\\
x_3
\\
\vdots 
\\
x_i
\end{array}
\right)
+
\left(
\begin{array}{c}
\theta_{ r +  ij + 1 }
\\
\theta_{ r + ij + 2 }
\\
\theta_{ r + ij + 3 }
\\
\vdots 
\\
\theta_{ r + ij + j }
\end{array}
\right),
\end{equation}
and ${\cal L}_l \colon \R^l \to \R^l$ is a mapping of the form 
\begin{equation}
{\cal L}_l(x_1, \dots, x_l) = (\rho(x_1), \dots, \rho(x_l))
\end{equation}
for a weakly differentiable continuous function $\rho \colon \R \to \R$.

\eqref{eq:Vn} is a feedforward neural network with activation function $\rho$ and 
$k+1$ layers (an input layer with $d$ neurons, $k-1$ hidden layers with $l$ neurons each, and an output 
layer with one neuron); see, e.g., \cite{Bengio09,LeCunBengioHinton15}.
Commonly used activation functions are e.g., the logistic function $x \mapsto e^x/(1+ e^x)$ or the 
ReLU function $x \mapsto \max \{0,x\}$.
The logistic function is continuously differentiable. So the corresponding neural networks $V_n$ have
well-defined $\theta$- and $x$-gradients $\nabla_{\theta} V_n$ and $\nabla_x V_n$. 
However, in the examples of Section \ref{sec:examples} below we use the ReLU function. It is 
continuously differentiable on $(-\infty,0) \cup (0, \infty)$ and has a left-hand derivative at $0$.
This yields weak $\theta$- and $x$-gradients for $V_n$, which we also denote by $\nabla_{\theta} V_n$ and $\nabla_x V_n$, 
respectively.

\subsection{Stochastic gradient descent based minimization}

Next, we recursively solve quadratic minimization problems 
to find parameter vectors $\theta^{1},\theta^{2},\ldots,$ $\theta^{N}\in\R^{\nu}$  
such that $V_n(\theta^n, x) \approx v(t_n, x)$ for $n \in \{1, 2, \dots, N\}$. More specifically,
$\theta^0 \in \R^{\nu}$ can be chosen arbitrarily; e.g., $\theta^0 = (0, \dots, 0) \in \R^{\nu}$. For
given $n\in\{1,2,\ldots,N\}$ and $\theta^0,\theta^1,\dots,\theta^{n-1} \in \R^\nu$,
one tries to find an approximate minimizer $\theta^{n}\in\R^{\nu}$ of the function 
\begin{multline}\label{eq:toMinimize}
%\begin{split}
% & 
\R^{\nu} \ni \theta
\mapsto 
\E\Big[\big|
V_{n}(\theta,Y_{T-t_{n}}) 
- 
\big[
V_{n-1}(\theta^{n-1},Y_{T-t_{n-1}}) 
\\
% & \quad 
+ f \big(Y_{T-t_{n-1}}, V_{n-1}(\theta^{n-1},Y_{T-t_{n-1}}),
\nabla_x V_{n-1}(\theta^{n-1},Y_{T-t_{n-1}})\big)\,(t_{n}-t_{n-1}) 
\big]
\big|^2\Big] \in \R.
%\end{split} 
\end{multline}
To do this with a standard stochastic gradient descent, one can initialize $\vartheta^n_0$ randomly or, e.g., 
as $\vartheta^n_0 = \theta^{n-1}$. Then one chooses a stepsize $\gamma \in (0,\infty)$ together with a number 
$M \in \N$ and iteratively updates for $m \in \{0,1, \dots, M-1\}$ according to
\begin{multline}\label{eq:sgdWithOriginalY}
%\begin{split}
\vartheta^{n}_{m+1} 
% & 
=  
\vartheta^{n}_m 
- 
2 \, \gamma \,
\nabla_{\theta} V_n(\vartheta^{n}_m,Y^{m}_{T-t_n})
\Big[ V_n(\vartheta^{n}_m,Y^{m}_{T-t_n}) 
-  V_{n-1}(\theta^{n-1},Y^{m}_{T-t_{n-1}}) \\
%& \quad 
-f \big(Y^{m}_{T-t_{n-1}},V_{n-1}(\theta^{n-1},Y^{m}_{T-t_{n-1}}),
\nabla_x V_{n-1}(\theta^{n-1}, Y^{m}_{T-t_{n-1}})\big)\,(t_{n}-t_{n-1})
\Big], 
% \end{split}
\end{multline}
where $Y^{m} \colon [0,T] \times \Omega \to \R^d$, $m \in \{0,1,\dots, M-1\}$, are 
$(\mathcal{F}_t)_{t\in [0,T]}$-adapted stochastic process with continuous sample paths
satisfying for every $t \in [0,T]$, the SDEs
% solving the SDEs
\begin{equation}
\label{eq:SDE-Y-m}
Y^{m}_t = \xi^{m} + \int_0^t \mu(Y^{m}_s)\,ds 
+ \int_0^t \sigma(Y^{m}_s)\,dB^{m}_s \quad \mbox{
	%for every $t \in [0,T]$ 
	$\P$-a.s.}
\end{equation}
corresponding to i.i.d.\ standard $(\mathcal{F}_t)_{t\in [0,T]}$-Brownian motions $B^{m}\colon [0,T]\times\Omega\to\R^d$
and i.i.d.\ $\mathcal{F}_0/\B(\R^d)$-measurable functions $\xi^{m}\colon\Omega\to\R^d$, $m \in \{0,1,\dots, M-1\}$.
After $M$ gradient steps one sets $\theta^n = \vartheta^n_M$.

\subsection{Discretization of the auxiliary stochastic process \texorpdfstring{$Y$}{Y}}
\label{subsec:Discrete-Y}

Equation \eqref{eq:sgdWithOriginalY} provides an implementable numerical algorithm in the special case where 
the solutions $Y^m$ of the SDEs \eqref{eq:SDE-Y-m} can be simulated exactly. If this is not the case, 
one can use a numerical approximation method to approximatively simulate $Y^m$, $m\in \{0,1,\dots, M-1\}$.
In the following we concentrate on the Euler--Maruyama scheme. But it is also possible to use a different 
approximation method. 

Note that it follows from \eqref{eq:SDE-Y-m} that
for  all $m\in \{0,1,\dots, M-1\}$ and $n\in \{0,1,\ldots,N-1\}$ we have
\begin{equation}
Y^{m}_{T-t_{n}} 
= 
Y^{m}_{T-t_{n+1}} 
+ 
\int_{T-t_{n+1}}^{T-t_n} \mu(Y^{m}_s)\,ds
+
\int_{T-t_{n+1}}^{T-t_n} \sigma(Y^{m}_s)\,dB^{m}_s \quad \mbox{$\P$-a.s.,}
\end{equation}
or equivalently, 
\begin{equation}
Y^m_{\tau_{n+1}} = Y^m_{\tau_{n}}  + \int_{\tau_n}^{\tau_{n+1}} \mu(Y^m_s)\,ds  + \int_{\tau_n}^{\tau_{n+1}} 
\sigma(Y^m_s)\,dB^m_s \quad \mbox{$\P$-a.s.}
\end{equation}
for $\tau_n = T - t_{N-n}$.
This suggests that for all $m \in \{0,1,\dots, M-1\}$ and $n \in \{0,1,\dots,N-1\}$, 
\begin{equation}
\label{NR3}
Y^m_{\tau_{n+1}} \approx Y^m_{\tau_{n}}  + \mu( Y^m_{\tau_{n}}) \,(\tau_{n+1}-\tau_n)
+ \sigma( Y^m_{\tau_{n}})\, (B^m_{\tau_{n+1}}-B^m_{\tau_n}).
\end{equation}
Therefore, we introduce the Euler--Maruyama approximations $\mathcal{Y}^m \colon \{0,1,\dots,N\} \times \Omega \to \R^d$ 
given for every $n \in \{0,1, \dots, N-1\}$ by $\Y^m_0 = \xi^m$ and %for every $n \in \{0,1, \dots, N-1\}$ that
\begin{equation}
\label{NR4}
\Y^{m}_{n+1} 
= 
\Y^{m}_{n} 
+ 
\mu(\Y^{m}_{n})\,(\tau_{n+1}-\tau_n)
+ 
\sigma(\Y^{m}_{n})\,(B^{m}_{\tau_{n+1}} - B^{m}_{\tau_n}).
%, \quad n \in \{0,1, \dots, N-1\}.
\end{equation}
It can be seen from \eqref{NR3} and \eqref{NR4} that for every $n \in \{0,1, \dots, N\}$,
one has 
\begin{equation} 
\Y^m_n\approx Y^m_{\tau_n}= Y^m_{T-{t_{N-n}}} \quad \mbox{and hence,} \quad
Y^m_{T-{t_{n}}} \approx \Y^m_{N-n},
% \quad n \in \{0,1, \dots, N\},
\end{equation}
which can be used to derive approximations of $(\vartheta^n_m)_{m=0}^M$,
$n \in \{1,2,\dots,N\}$, from \eqref{eq:sgdWithOriginalY} that are also implementable if the 
processes $Y^m$ cannot be simulated exactly. More precisely, set
$\Theta^0_M = (0, \dots, 0) \in \R^{\nu}$. For $n \in \{1,2,\dots, N\}$, initialize $\Theta^n_0$
randomly or as $\Theta^n_0 = \Theta^{n-1}_M$. Then set for every  $m \in \{0, 1, \dots, M-1\}$,
\begin{multline}\label{eq:sgdWithApproximatedY}
%\begin{split}
\Theta^{n}_{m+1} = 
\Theta^{n}_m 
- 2 \,\gamma \,
\nabla_{\theta} V_n(\Theta^{n}_m,\Y^{m}_{N-n})
\Big[
V_n(\Theta^{n}_m,\Y^{m}_{N-n}) 
- V_{n-1}(\Theta^{n-1}_M,\Y^{m}_{N-n+1})\\ 
%  & \quad 
-f\big(\Y^{m}_{N-n+1}, V_{n-1}(\Theta^{n-1}_M,\Y^{m}_{N-n+1}),
\nabla_x V_{n-1}(\Theta^{n-1}_M,\Y^{m}_{N-n+1})\big)\,(t_{n}-t_{n-1})
\Big].
% \end{split}
\end{multline}
Comparing \eqref{eq:sgdWithApproximatedY} to \eqref{eq:sgdWithOriginalY} suggests that $\Theta^n_m \approx \vartheta^n_m$ 
for $m \in \{0,1, \dots, M\}$ and $n \in \{1, 2, \dots, N\}$.

In the following two Subsections \ref{subsec:algo1} and \ref{subsec:algo-Full-gen},
we first formalize and then generalize the approximation algorithm
derived in Subsections \eqref{subsec:temp-discret}--\eqref{subsec:Discrete-Y}.

\subsection{Description of the algorithm in a special case}
\label{subsec:algo1}

In this subsection, we provide a formal description of the algorithm 
derived in Subsections \eqref{subsec:temp-discret}--\eqref{subsec:Discrete-Y}
in the case where the standard Euler--Maruyama scheme (cf., e.g., \cite{KloedenPlaten1992,Maruyama1955,MilsteinOriginal1974}) 
is used to approximate the solutions $Y^m$ of the SDEs \eqref{eq:SDE-Y-m} and
optimal parameters $\theta_1, \dots, \theta_N \in \R^{\nu}$ are computed with 
plain vanilla stochastic gradient descent with a constant learning rate $\gamma \in (0,\infty)$ 
and batch size 1. Note that for the algorithm to be implementable, it is enough 
if the initial condition $\varphi$ has a weak gradient $\nabla \varphi \colon \R^d \to \R^d$
which does not need to satisfy any growth conditions.

In the following Framework~\ref{algo:1}, feedforward neural networks of the form \eqref{eq:Vn} are employed
to approximate the solution of the PDE.  \eqref{eq:plain-vanilla-SGD} describes a 
stochastic gradient descent scheme with constant learning rate $\gamma$ and
\eqref{eq:algo-1-quadratic-miniminzation} specifies the quadratic loss functions.
\begin{algo}[Special case]
\label{algo:1}
Assume $\varphi$ has a weak gradient $\nabla \varphi \colon \R^d \to \R^d$.
Consider $N \in \N$ and $t_0,t_1,\ldots,t_N\in [0,T]$ such that 
\begin{equation}\label{eq:algo1-time}
0 = t_0 < t_1 < \ldots < t_N = T.
\end{equation}
Set $\tau_n = T-t_{N-n}$ for $n \in \{0,1,\dots, N\}$. For a given $M \in \N$, let $B^{m}\colon [0,T]\times\Omega\to\R^d$, 
$m\in \{0,1,\dots, M-1\}$, be i.i.d.\ standard Brownian motions on a probability space $(\Omega,\F, \P)$.
Consider i.i.d.\ random variables $\xi^m \colon \Omega \to \R^d$, $m \in \{0,1,\dots, M-1\}$,
that are independent of $B^m$, $m \in \{0,1,\dots, M-1\}$,
and let the stochastic processes $\Y^{m} \colon \{0,1,\ldots,N\} \times\Omega\to\R^d$, $m \in \{0,1,\dots, M-1\}$, 
be 
given by $\Y^{m}_0 = \xi^{m}$ and
\begin{equation}\label{Y-algo-spez}
	\Y^{m}_{n+1} 
	= 
	\Y^{m}_n 
	+ 
	\mu(\Y^{m}_n)\,(\tau_{n+1}-\tau_{n})
	+ 
	\sigma(\Y^{m}_n)\,(B^{m}_{\tau_{n+1}}-B^{m}_{\tau_{n}}), \quad n \in \{0, 1, \dots, N-1\}.
	\end{equation}  
Let $V_n \colon \R^{\nu} \times \R^d \to \R$, $n \in \{0,1,\dots, N\}$ be the functions given in 
\eqref{eq:V0}--\eqref{eq:Vn}. Consider $\gamma \in (0,\infty)$ and let
$\Theta^{n}\colon \{0, 1, \dots, M\} \times\Omega\to\R^{\nu}$, $n \in \{0,1,\ldots,N\}$, 
be  stochastic processes satisfying for all $m \in \{0,1,\dots, M-1\}$ and $n \in \{1, \dots, N\}$,
\begin{equation}
\label{eq:plain-vanilla-SGD}
\Theta^{n}_{m+1} = \Theta^{n}_m - \gamma \, \Phi^{n,m}(\Theta^{n}_m),
\end{equation}
where for all $\omega \in \Omega$,
\begin{equation}
		\label{eq:algo-1-quadratic-miniminzation}
	\begin{split}
	%    & 
	&\phi^{n,m}(\theta,\omega) 
	= \Big[ V_n\big(\theta,\Y^{m}_{N-n}(\omega)\big) - V_{n-1}(\Theta^{n-1}_M(\omega),\Y^{m}_{N-n+1}(\omega)) 
	\, - (t_{n}-t_{n-1}) \\
	& \times f\big(\Y^{m}_{N-n+1}(\omega), V_{n-1}(\Theta^{n-1}_M(\omega),\Y^{m}_{N-n+1}(\omega)),
	\nabla_x V_{n-1}(\Theta^{n-1}_M(\omega),\Y^{m}_{N-n+1}(\omega))\big) \Big]^2, %\nonumber
	\end{split}
	\end{equation}
	and
	$\Phi^{n,m}(\theta,\omega) = \nabla_{\theta}\phi^{n,m}(\theta,\omega)$.\end{algo}
In the setting of Framework~\ref{algo:1} the random variables $V_n(\Theta^n_M,x) \colon \Omega \to \R$
provide for all $n \in \{1,2, \dots, N\}$ and $x \in \R^d$ the approximations 
\begin{equation}
V_n(\Theta^n_M,x) \approx u(t_n,x)
\end{equation}
%$n \in \{1,2, \dots, N\}$, $x \in \R^d$, 
to the solution $u \colon [0,T] \times \R^d \to \R$ of the PDE 
\eqref{eq:defPDE}.

\subsection{Description of the algorithm in the general case}
\label{subsec:algo-Full-gen}
%In this subsection we provide a general framework (see Framework~\ref{def:general_algorithm} below) which covers the deep splitting method derived in Subsections~\ref{subsec:temp-discret}--\ref{subsec:algo1}. The variant of the deep splitting method described in Subsection~\ref{subsec:algo1} still remains the core idea of Framework~\ref{def:general_algorithm}. However, Framework~\ref{def:general_algorithm} allows to incorporate other minimization algorithms 
%
We now generalize Framework~\ref{algo:1} so that, besides plain vanilla stochastic gradient descent
with constant learning rate and batch size $1$, it also covers more advanced machine learning techniques
such as mini-batches, batch normalization
and more sophisticated updating rules.

In Framework~\ref{algo2} below, functions $V^{j,s}_n\colon\R^\nu\times\R^d \to\R$ 
parametrized by $(j,s,n) \in \N \times \R^\varsigma \times\{0,1,\dots,N\}$ 
for some number $\varsigma \in \N$, are used to approximate the solution of the PDE.
The additional parameters $j$ and $s$ make it possible 
to describe mini-batches and batch normalization; see Ioffe \& Szegedy \cite{IoffeSzegedy2015}.
As in Framework~\ref{algo:1}, the standard Euler--Maruyama scheme 
is used to approximate the solutions $Y^m$ of the SDEs \eqref{eq:SDE-Y-m}.
But a different approximation scheme could be employed as well.
In \eqref{def:general_algorithm-quadratic-loss} the quadratic loss functions are given
that are used for training the functions $V^{j,s}_n$, whereas \eqref{def:general_algorithm-quadratic-loss-grad}
specifies the gradients of the loss functions.
The role of the stochastic processes $\mathbb{S}^n$ in \eqref{eq:general_batch_normalization} is to 
describe the variables (running mean and standard deviation) needed for batch normalization.
The stochastic processes $\Theta^{n}$ and $\Xi^{n}$ in \eqref{eq:general_gradient_step} describe the 
updating rule. Their dynamics are specified by the functions $\Psi^n_m$ and $\psi^n_m$.
Since we make no assumptions on $\Psi^n_m$ and $\psi^n_m$, the framework includes 
very general stochastic gradient optimization methods. In the examples in Section \ref{sec:examples}
we use Adam optimization; see Kingma \& Ba~\cite{KingmaBa2015}. The corresponding 
specification of the functions $\Psi^n_m$ and $\psi^n_m$ is given in 
\eqref{eq:examples_setting_moment_estimation}--\eqref{eq:examples_setting_adam_grad_update} below.
\begin{algo}[General case] 
\label{algo2}
Assume $\varphi$ has a weak gradient $\nabla \varphi \colon \R^d \to \R^d$. Consider $N \in \N$ and 
$t_0,t_1,\dots,t_N \in [0,T]$ such that $0 = t_0 < t_1 < \ldots < t_N = T$. Set $\tau_n= T-t_{N-n}$ for 
$n \in \{0,1,\dots,N\}$. For a given $M \in \N$, let, for every $n \in \{1,2,\ldots,N\}$,
$B^{n,m,j} \colon [0,T] \times \Omega \to \R^d$, 
$m \in \{0,1, \dots, M-1\}$, $j \in \N$,
be i.i.d.\ standard Brownian motions on a probability space $(\Omega,\F,\P)$
and $\xi^{n,m,j} \colon \Omega\to\R^d$, $m \in \{0,1, \dots, M-1\}$, $j \in \N$, i.i.d.\ random variables
that are independent of $B^{n,m,j}$, $m \in \{0,1, \dots, M-1\}$, $j \in \N$.
Let $\Y^{n,m,j}\colon \{0,1,\ldots,N\}\times\Omega\to\R^d$, $n \in \{1,2,\ldots,N\}$, 
$ m \in \{0,1,\dots, M-1\} $, $ j \in \N $, be stochastic processes given by $\Y^{n,m,j}_0 = \xi^{n,m,j}$ and 
\begin{equation}\label{eq:FormalXapprox}
	\Y^{n,m,j}_{k+1} 
	= \mu(\Y^{n,m,j}_k)\,(\tau_{k+1}-\tau_{k}) + 
		\sigma(\Y^{n,m,j}_k)\,(B^{n,m,j}_{\tau_{k+1}}-B^{n,m,j}_{\tau_{k}}), \quad
	k \in \{0, 1, \dots, N-1\}.
		\end{equation}  
Let $\nu, \varsigma, \varrho, J_0, \dots, J_{M-1} \in \N$. Consider functions $V^{j,s}_n \colon\R^\nu\times\R^d \to\R$, 
$(j,s,n) \in \N \times \R^\varsigma \times\{0,1,\dots,N\}$, such that $V^{j,s}_0(\theta,x) = \varphi(x)$ 
for all $(j,s,\theta,x) \in \N \times \R^{\varsigma} \times \R^{\nu} \times \R^d$, and let
$\Theta^{n}\colon \{0,1, \dots, M-1\} \times\Omega\to\R^{\nu}$, $n \in \{0,1,\dots N\}$, be stochastic processes.
For all  $n\in\{1,2,\ldots,N\}$, $m\in \{0,1, \dots, M-1\}$ and $s \in\R^{\varsigma}$, let the mapping 
$\phi^{n,m,s}\colon\R^{\nu}\times\Omega\to\R$ be given by 
	\begin{multline}
		\label{def:general_algorithm-quadratic-loss}
	%\begin{split}
	\phi^{n,m,s}(\theta,\omega) 
	= 
	\frac{1}{J_m}\sum_{j=1}^{J_m}
	\bigg[ 
	V^{j,s}_n\big(\theta,\Y^{n,m,j}_{N-n}(\omega)\big) 
	- 
	V^{j,s}_{n-1}\big(\Theta^{n-1}_{M}(\omega),\Y^{n,m,j}_{N-n+1}(\omega)\bigr) - (t_n-t_{n-1})  
	\\
	%& \quad \quad \cdot
	\times f \Big(\Y^{n,m,j}_{N-n+1}(\omega),
	V^{j,s}_{n-1}\big(\Theta^{n-1}_{M}(\omega),\Y^{n,m,j}_{N-n+1}(\omega)\bigr), 
	\nabla_x V^{j,s}_{n-1}\big(\Theta^{n-1}_{M}(\omega),\Y^{n,m,j}_{N-n+1}(\omega)\big)\Big)
	\bigg]^2,
	\end{multline}
and assume $\Phi^{n,m,s }\colon\R^{\nu}\times\Omega\to\R^{\nu}$ is a function satisfying 
\begin{align}
		\label{def:general_algorithm-quadratic-loss-grad}
	\Phi^{n,m,s}(\theta,\omega) = \nabla_{\theta}\phi^{n,m,s}(\theta,\omega)
	\end{align}
for all $\omega\in\Omega$ and $\theta \in\{\vartheta\in\R^{\nu}\colon \phi^{n,m,s}(\cdot,\omega)\colon\R^{\nu}\to\R~\text{is differentiable at}~\vartheta\}$. For every $n\in\{1,2,\ldots,N\}$ and $m \in \{0,1,\dots, M-1\}$, let 
$\S^n\colon\R^{\varsigma}\times\R^{\nu}\times(\R^d)^{\{0,1,\ldots,N\}\times\N}\to\R^{\varsigma}$,
$\Psi^n_m\colon\R^{\varrho}\times\R^{\nu}\to\R^{\varrho}$ 
and $\psi^n_m\colon\R^{\varrho}\to\R^{\nu}$ be functions and
$\bS^{n}\colon \{0,1,\dots, M-1\} \times\Omega\to\R^{\varsigma}$ 
and  $\Xi^{n}\colon \{0,1,\dots, M-1\} \times\Omega\to\R^{\varrho}$ stochastic processes such that
	\begin{equation}\label{eq:general_batch_normalization} 
	\bS^{n}_{m+1} = \S^{n}\bigl(\bS^{n}_m, \Theta^{n}_{m}, 
	(\Y_k^{n,m,i})_{(k,i)\in\{0,1,\ldots,N\}\times\N}\bigr),  
	\end{equation}
	\begin{equation}
	\Xi^n_{m+1} = \Psi^n_{m}(\Xi^n_{m},\Phi^{n,m,\bS^n_{m+1}}(\Theta^n_m))
	\quad
	\text{and}
	\quad
	\Theta^{n}_{m+1} = \Theta^{n}_{m} - \psi^n_{m}(\Xi^n_{m+1}) 
	\label{eq:general_gradient_step}. 
	\end{equation}
\end{algo}
In the setting of Framework~\ref{algo2} the functions 
$V^{1,\mathbb{S}_M^n}_n(\Theta^n_M,x)\colon \Omega \to \R$ yield the approximations 
\begin{equation}\label{eq:Algo2-PDE-approx}
V^{1,\mathbb{S}_M^n}_n(\Theta^n_M,x) \approx u(t_n,x), \quad n \in \{1, 2, \dots, N\}, \; x \in \R^d,
\end{equation}
%$n \in \{1, 2, \dots, N\}$, $x \in \R^d$,
of the solution $u \colon [0,T] \times \R^d \to $ to the PDE \eqref{eq:defPDE}.
%
%
%%%%%%%%%%%%%%%%
%The role of the processes $\mathbb{S}^n\colon \N_0\times\Omega \to \R^{\varsigma}$, $n\in\{1,2,\ldots,N\}$, is to describe the variables needed for batch normalization.
%%%%%%%%%%%%%%%%%%%%%%%%%%%%%%%%%%%%%%%%%%%%%%%%%%%%%%%%%%%%
%
%
\section{Examples}
\label{sec:examples}
%\input{Examples.tex}
%
%%%%%%%%%%%%%%%%%%%%%%%%%%
We now illustrate the performance of the deep splitting method on five concrete example PDEs. 
In each example we use the general approximation method of Framework \ref{algo2}
with approximating functions $V^{j,s}_n \colon \R^{\nu} \times \R^d \to \R$, $n \in \{1,2, \dots, N\}$,
specified as feedforward neural networks with 4 layers (1 input layer, 2 hidden layers, 1 output layer)
and ReLU-activation $\rho(x) = \max\{0,x\}$, $x \in \R$. 
We use mini-batches of size $J_m = 256$ and apply batch normalization before the
first affine transformation, before each of the two nonlinear activation functions in front of the hidden layers, 
and just before the output layer. We use Xavier initialization (see Glorot \& Bengio \cite{glorot2010understanding}) 
to initialize all weights in the neural networks together with Adam optimization (see Kingma \& Ba~\cite{KingmaBa2015})
with parameters $\varepsilon=10^{-8}$, $\beta_1 = 0.9$, $\beta_2 = 0.999$ and 
decreasing learning rates $(\gamma_m)_{m =0}^{M-1}$ that we choose depending on the form 
and dimension of the problem. More precisely, we set $\varrho = 2\nu$ and denote by 
$\operatorname{Pow}_2 \colon \R^{\nu} \to\R^{\nu}$ the function given by 
$\operatorname{Pow}_2(\eta_1, \dots, \eta_{\nu}) =  (\eta_1^2, \dots, \eta^2_{\nu})$. Then,
Adam optimization corresponds to the following specification of the 
two functions $\Psi^n_m\colon \R^{3\nu}\to\R^{2\nu}$ 
and $\psi^n_m \colon\R^{2\nu}\to\R^{\nu}$ from Framework \eqref{algo2}:

\begin{align}\label{eq:examples_setting_moment_estimation}
	\Psi^n_m ( x , y , \eta ) 
	= 
	(\beta_1 x + (1-\beta_1) \eta \, , \, \beta_2 y + (1-\beta_2) \operatorname{Pow}_2(\eta)) %,
	\end{align}
	and 
	\begin{align}\label{eq:examples_setting_adam_grad_update}
	\psi^n_m ( x,y ) = 
	\biggl(
	\Bigl[
	\sqrt{\tfrac{|y_1|}{1- \beta_2^m}} + \varepsilon
	\Bigr]^{-1}
	\frac{\gamma_m x_{1}}{1- \beta_1^m},
	\ldots, 
	\Bigl[
	\sqrt{\tfrac{|y_{\nu}|}{1- \beta_2^m}} + \varepsilon
	\Bigr]^{-1}
	\frac{\gamma_m x_{\nu}}{1- \beta_1^m}
	\biggr). 
	\end{align}

In the examples in the following subsections, we approximate $u(T,x)$ 
for different $T \in (0,\infty)$ and $x \in \R^d$ with $V^{1,\mathbb{S}^N_M}_N(\Theta^N_M,x)$,
which, due to the stochastic gradient optimization method, is a random variable. In each example we report
estimates of the expectation and standard deviation of $V^{1,\mathbb{S}^N_M}_N(\Theta^N_M,x)$. 
We also give relative $L^1$-approximation errors with respect to reference values
calculated with different alternative methods together with their uncorrected sample standard deviations.
The average runtimes needed for calculating one realization of $V^{1,\mathbb{S}^N_M}_N(\Theta^N_M,x)$
are in seconds and were determined as averages over 10 independent runs.

All numerical experiments presented in this paper were implemented in {\sc Python} using {\sc TensorFlow}
and run on a NVIDIA GeForce GTX 1080 GPU with 1974 MHz core clock and 8 GB
GDDR5X memory with 1809.5 MHz clock rate and an underlying system consisting of
an Intel Core i7-6800K 3.4 GHz CPU with 64 GB DDR4-2133 memory running TensorFlow 1.5 on Ubuntu 16.04. 
The {\sc Python} source codes can be found at https://github.com/seb-becker/deep$\_$pde.
%
%We would like to point out that no special emphasis has been
%put on optimizing computation speed. In many cases some of the algorithm parameters
%could be adjusted in order to obtain similarly accurate results in shorter runtime.
%%%%%%%%%%%%%%%%%%%%
%We also refer to 
%the {\sc Python} codes %\ref{code:deepPDEmethod} 
%in Section~\ref{sec:sourcecodes} below.
%%%%%%%%%%%%%%%%%%%%%%%%%%%%%%%%%%%%
% for an  
%implementation of the deep splitting method in the case of .

%
%%%%%%%%%%%%%%%%%%%%%%%%%%%%%%%%%%%%%%%%%%%%%%%%%%%%%%%%%%%%%%%%%%%%%%%%%%%%%%
%%%%%%%%%%%%%%%%%%%%%%%%%%%%%%%%%%%%%%%%%%%%%%%%%%%%%%%%%%%%%%%%%%%%%%%
\subsection{Hamilton--Jacobi--Bellman (HJB) equations}
\label{subsec:HJB}
%%%%%%%%%%%%%%%%%%%%%%%%%%%%%%%%%%%%%%%%%%%%%%%%%%%%%%%
%\input{table_HamiltonJacobiBellman.tex}
%%%%%%%%%%%%%%%%%%%%%%%%%%%%%%%
%%%%%%%%%%%%%%%%%%%%%%%%%%%%%%%%%%%%%%%%%%%%%%%%%%%%%% 
\begin{table}
	\begin{center}
		\resizebox{\textwidth}{!}{\begin{tabular}{|c|c|c|c|c|c|c|c|c|}
				\hline
				$ d $ & $ T $ & $ N $ & Expectation & Std.\ dev.\ & Ref.\ value & rel.\ $L^1$-error & Std.\ dev.\ rel.\ error & avg.\ runtime\\
				\hline
				10  &  \nicefrac{1}{3} & 8 & 1.56645 & 0.00246699 & 1.56006 & 0.00410 & 0.00158134 & 18.0s \\
				10  &  \nicefrac{2}{3} & 16 & 1.86402 & 0.00338646 & 1.85150 & 0.00677 & 0.00182904 & 37.9s \\
				10  &  1 & 24 & 2.07017 & 0.00634850 & 2.04629 & 0.01167 & 0.00310245 & 58.2s \\
				\hline
				50  &  \nicefrac{1}{3} & 8 & 2.39214 & 0.00151918 & 2.38654 & 0.00234 & 0.00063656 & 18.0s \\
				50  &  \nicefrac{2}{3} & 16 & 2.84607 & 0.00140300 & 2.83647 & 0.00338 & 0.00049463 & 37.9s \\
				50  &  1 & 24 & 3.15098 & 0.00275839 & 3.13788 & 0.00417 & 0.00087906 & 58.4s \\
				\hline
				100  &  \nicefrac{1}{3} & 8 & 2.85090 & 0.00071267 & 2.84696 & 0.00138 & 0.00025033 & 18.1s \\
				100  &  \nicefrac{2}{3} & 16 & 3.39109 & 0.00093368 & 3.38450 & 0.00195 & 0.00027587 & 38.2s \\
				100  &  1 & 24 & 3.75329 & 0.00136920 & 3.74471 & 0.00229 & 0.00036564 & 58.3s \\
				\hline
				200  &  \nicefrac{1}{3} & 8 & 3.39423 & 0.00051028 & 3.39129 & 0.00087 & 0.00015047 & 18.1s \\
				200  &  \nicefrac{2}{3} & 16 & 4.03680 & 0.00088215 & 4.03217 & 0.00115 & 0.00021878 & 38.0s \\
				200  &  1 & 24 & 4.46734 & 0.00079688 & 4.46172 & 0.00126 & 0.00017860 & 58.2s \\
				\hline
				300  &  \nicefrac{1}{3} & 8 & 3.75741 & 0.00063334 & 3.75530 & 0.00056 & 0.00016865 & 18.3s \\
				300  &  \nicefrac{2}{3} & 16 & 4.46859 & 0.00049953 & 4.46514 & 0.00077 & 0.00011187 & 38.5s \\
				300  &  1 & 24 & 4.94586 & 0.00087736 & 4.94105 & 0.00097 & 0.00017756 & 58.8s \\
				\hline
				500  &  \nicefrac{1}{3} & 8 & 4.27079 & 0.00051256 & 4.26900 & 0.00042 & 0.00012007 & 18.0s \\
				500  &  \nicefrac{2}{3} & 16 & 5.07900 & 0.00034792 & 5.07618 & 0.00056 & 0.00006854 & 38.0s \\
				500 &  1 & 24 & 5.62126 & 0.00045092 & 5.61735 & 0.00070 & 0.00008027 & 57.7s \\
				\hline
				1000  &  \nicefrac{1}{3} & 8 & 5.07989 & 0.00022764 & 5.07876 & 0.00022 & 0.00004482 & 20.6s \\
				1000  &  \nicefrac{2}{3} & 16 & 6.04130 & 0.00030680 & 6.03933 & 0.00033 & 0.00005080 & 43.7s \\
				1000  &  1 & 24 & 6.68594 & 0.00040334 & 6.68335 & 0.00039 & 0.00006035 & 66.5s \\
				\hline
				5,000  &  \nicefrac{1}{3} & 8 & 7.59772 & 0.00024745 & 7.59733 & 0.00005 & 0.00003257 & 120.4s \\
				5,000  &  \nicefrac{2}{3} & 16 & 9.03721 & 0.00027322 & 9.03466 & 0.00028 & 0.00003024 & 256.7s \\
				5,000  &  1 & 24 & 9.97266 & 0.00047098 & 9.99835 & 0.00257 & 0.00004711 & 393.9s \\
				\hline
				10,000  &  \nicefrac{1}{3} & 8 & 9.03574 & 0.00022994 & 9.03535 & 0.00004 & 0.00002545 & 519.8s \\
				10,000  &  \nicefrac{2}{3} & 16 & 10.74521 & 0.00026228 & 10.74478 & 0.00004 & 0.00002157 & 1105.6s \\
				10,000  &  1 & 24 & 11.87860 & 0.00022705 & 11.89099 & 0.00104 & 0.00001909 & 1687.7s \\
				\hline
		\end{tabular}}
		\caption{Deep splitting approximations of the solution of the HJB equation \eqref{eq:HJB}
		for different values of $d$, $T$ and $N$.}		\label{table:HJB}
	\end{center}
\end{table}
%%%%%%%%%%%%%%%%%%%%%%%%%%%%%%%%%%%%%%%%%%%%%%%%%%%%%%%%%%%%
%
% 
%\begin{figure}
%	\includegraphics{HamiltonJacobiBellman.pdf}
%%%%%%%%%%%%%%%%%%%%%%%%%
%	\caption{Plots of approximative 
%		calculations of the relative $L^1$-approximation error 
%		$\E\big[\tfrac{|\U^{\Theta_m}(\xi)-4.5901|}{4.5901}\big]$ 
%		and of the mean of the empirical loss function $\E\big[\tfrac{1}{J_m}
%		\sum_{j=1}^{J_m}|\Y^{m,\Theta_m,j,\bS_{m+1}}_N-g(\X^{m,j}_N)|^2\big]$ 
%		in the case of the $100$-dimensional 
%		Hamilton-Jacobi-Bellman equation \eqref{eq:example_hamilton_jacobi_bellman} 
%		against $ m\in\{0,1,\ldots,2000\} $.} 
%	\label{fig:figure_HamiltonJacobiBellman}
%	%%%%%%%%%%%%%%%%%%%%%%%%%%%%%%%%%%%%%%%%%%%%%%%%%%
%\end{figure}
%%%%%%%%%%%%%%%%%%%%%%%%%%%%%%%%%%%%%%%%%%%%%%%%%%%%%%
As a first example, we calculate approximate solutions of the PDE %Framework~\ref{def:general_algorithm}
\begin{equation}\label{eq:HJB}
 \tfrac{\partial}{\partial t} u(t,x) = 
 \Delta_x u(t,x) - \|\nabla_x u(t,x)\|_{\R^d}^2, \quad (t, x) \in (0,T] \times \R^d,
\end{equation} 
with initial condition $u(0,x) = \|x\|^{1/2}_{\R^d}$ for different $d\in\N$. 
The deep splitting method can be applied to more general HJB equations.
But \eqref{eq:HJB} has the advantage that it reduces to a linear heat equation under a logarithmic transformation;
see, e.g., E et al.~\cite[Lemma 4.2]{EHanJentzen2017}. Since solutions of linear equations can 
be approximated with standard Monte Carlo, this allows us to efficiently compute reference 
solutions in high dimensions. 

Table \ref{table:HJB} shows deep splitting approximations of $u(T,0, \dots, 0)$ for different values of 
$d$, $T$ and $N$. We derived it with $\mu(x)=(0,0,\dots,0) \in \R^d$, $\sigma(x)=\sqrt{2}\mbox{ Id}_{\R^{d\times d}}$, 
and $f(x,y,z)= -\|z\|^2_{\R^d}$ for $(x,y,z) \in \R^d \times \R \times \R^d$. We used
$M= 500 +  100\,\mathbbm{1}_{\{10,000\}}(d)$ and
\begin{equation}
\gamma_m = \begin{cases}10^{-1} \mathbbm{1}_{[0,300]}(m) 
+ 10^{-2} \mathbbm{1}_{(300,400]}(m)
+ 10^{-3} \mathbbm{1}_{(400,500]}(m) & \mbox{ for } d< \mbox{10,000}\\
10^{-1} \mathbbm{1}_{[0,400]}(m) 
+ 10^{-2} \mathbbm{1}_{(400,500]}(m)
+ 10^{-3} \mathbbm{1}_{(500,600]}(m) & \mbox{ for } d= \mbox{10,000}.
\end{cases}
\end{equation}
We set $\xi^{n,m,j}=(0,0,\dots,0) \in \R^{d}$ 
for every $(n,m,j) \in \{1,2, \dots, N\} \times \{0, 1, \dots, M-1\} \times \N$, and 
used feedforward neural networks with a $d$-dimensional 
input layer, two hidden layers of dimension $d+10$, and a one-dimensional output layer.

The reference values for $u(T,0,0,\dots,0)$ were calculated using a logarithmic transformation and 
a standard Monte Carlo method; see Han et al.\ \cite[Lemma~4.2]{EHanJentzen2017}.

%(see {\sc Matlab} code~ %\ref{code:branchingMatlab}
%in Appendix~ %\ref{subsec:BranchingMatlab} 
%below).
% 
% 
% 
%%%%%%%%%%%%%%%%%%%%%%%%%%%%%%%%%%%%%%%%%%%%%%%%%%%%%%%%%%%%%%%%%%%%%%%%%%%%%
%%%%%%%%%%%%%%%%%%%%%%%%%%%%%%%%%%%%%%%%%%%%%%%%%%%%%%%%%%%%%%%%%%%%%%

\subsection{Nonlinear Black--Scholes equations}
\label{subsec:NonLinBlackScholes}
%
%%%%%%%%%%%%%%%%%%%%%%%%%%%%%%%%%%%%%%%%%%%%%%%%%
%                                          \input{table_HamiltonJacobiBellman.tex}
%%%%%%%%%%%%%%%%%%%%%%%%%%%%%%%%%%%%%%%%%%%%%%%%%%%%%%
\begin{table}
	\begin{center}
		\resizebox{\textwidth}{!}{\begin{tabular}{|c|c|c|c|c|c|c|}
			\hline
			$ d $ & Expectation & Stdev & Ref. value & rel. $L^1$-error & Stdev rel. error & avg. runtime\\
			\hline
			10     & 40.6553107 & 0.1000347132 & 40.7611353 & 0.0029624273 & 0.0019393471 & 858.3s \\
			50     & 37.421057 & 0.0339765334 & 37.5217732 & 0.0026842068 & 0.0009055151 & 975.4s \\
			100    & 36.3498646 & 0.027989905 & 36.4084035 & 0.0016078403 & 0.000768776 & 1481.5s \\
			\hline
			200    & 35.374638 & 0.035236816 & 35.4127342 & 0.0012857962 & 0.0006625744 & 951.2s \\
			300    & 34.8476466 & 0.0225350305 & 34.8747946 & 0.0008818254 & 0.0004762554 & 953.3s \\
			500    & 34.2206181 & 0.0081072294 & 34.2357988 & 0.0004701552 & 0.0001701012 & 956.0s \\
			1000   & 33.4058827 & 0.0050161752 & 33.4358163 & 0.0008952555 & 0.000150024 & 1039.6s \\
			5,000   & 31.7511529 & 0.0048508218 & 31.7906594 & 0.0012427078 & 0.0001525864 & 7229.7s \\
			10,000  & 31.1215014 & 0.0031131196 & 31.1569116 & 0.0011365119 & 0.00009991746 & 23,593.2s \\
			\hline
		\end{tabular}}
		\caption{Deep splitting approximations of the solution of the nonlinear Black--Scholes equations 
			\eqref{eq:BS} for $T = 1/3$, $N= 96$ and different $d$.}
		\label{table:BS}
	\end{center}
\end{table}
%%%%%%%%%%%%%%%%%%%%%%%%%%%%%%%%%%%%%%%%%%%%%%%%%%%%
% 
% 
%\begin{figure}
%	\includegraphics{HamiltonJacobiBellman.pdf}
%%%%%%%%%%%%%%%%%%%%%%%%%%%%%%%%%%%%%%%%%%%%
%	\caption{Plots of approximative 
%		calculations of the relative $L^1$-approximation error 
%		$\E\big[\tfrac{|\U^{\Theta_m}(\xi)-4.5901|}{4.5901}\big]$ 
%		and of the mean of the empirical loss function $\E\big[\tfrac{1}{J_m}
%		\sum_{j=1}^{J_m}|\Y^{m,\Theta_m,j,\bS_{m+1}}_N-g(\X^{m,j}_N)|^2\big]$ 
%		in the case of the $100$-dimensional 
%		Hamilton-Jacobi-Bellman equation \eqref{eq:example_hamilton_jacobi_bellman} 
%		against $ m\in\{0,1,\ldots,2000\} $.} 
%	\label{fig:figure_HamiltonJacobiBellman}
%%%%%%%%%%%%%%%%%%%%%%%%%%%%%%%%%%%%%%%%%%%%%%%%%%%%%%%%%%%%%%%
%\end{figure}
%%%%%%%%%%%%%%%%%%%%%%%%%%%%%%%%%%%%%%%%%%%%%%%%%%%%%%%%
There exist a number of extensions of the classical linear Black--Scholes equation which incorporate
nonlinear phenomena such as transaction costs, default risk or Knightian uncertainty. We 
here consider nonlinear Black--Scholes equations of the form
\begin{equation}\label{eq:BS}
\begin{split}
\tfrac{\partial}{\partial t}u(t,x) &=  
- u(t,x)\,
(1-\delta)\left[\min\!\left\{
\gamma^h,\max\!\left\{\gamma^l,\tfrac{ ( \gamma^h - \gamma^l ) }{ ( v^h - v^l ) }
\left(u(t,x) - v^h \right)+\gamma^h\right\}
\right\}\right]
%  +
%  \int_{ 0 }^{ t }
%  b\big( x, X_s( x ), ( \nabla X_s )( x ) \big)
%  \, dZ_s(x)
\\
& \quad
- R u(t,x) + \big\langle \bar \mu\, x, \nabla_x u ( t,x ) \big\rangle_{ \R^d }
+
\tfrac{ \bar \sigma^2 }{ 2 }\!\left[ \textstyle\sum\limits_{i=1}^d |x_i|^2 \, \frac{\partial^2}{\partial x^2_i} u(t,x)\right], 
\end{split}
\end{equation}
$(t,x) \in (0,T] \times \R^d$, for suitable parameters $\delta$, $R$, $\gamma^h$, $\gamma^l$, 
$v^h$, $v^l$, $\bar{\mu}$, $\bar{\sigma} \in \R$. They describe derivative prices under 
default risk; see, e.g., Han et al.\ \cite[Subsection~3.1]{HanEJentzen17} and E et al.\ \cite[Subsection~3.1]{MultilevelPicard},
from where we adopt the initial condition $u(0,x)=\min_{i\in\{1,2,\dots,d\}} x_i$
and the choice of the parameter values $\delta = 2/3$, $R= 0.02$, $\gamma^h=0.2$, $\gamma^l=0.02$,
$v^h=50$, $v^l=70$, $\bar\mu=0.02$, $\bar \sigma=0.2$.

Table \ref{table:BS} reports deep splitting 
approximations of $u(T, 50, \dots, 50)$ for $T = 1/3$, $N = 96$, and different values of $d$.
We chose $\mu(x)= \bar \mu x $, $\sigma(x)=\bar \sigma x$, and 
\begin{equation} \label{eq:nbs_f}
f(x,y,z)= -(1-\delta)\min\!\left\{
\gamma^h,\max\left\{\gamma^l,\tfrac{ ( \gamma^h - \gamma^l ) }{ ( v^h - v^l ) }
\left(y - v^h \right)+\gamma^h\right\}
\right\}y-Ry,
\end{equation}
$(x,y,z) \in \R^d \times \R \times \R^d$ and used $M= 2000 +  1000\,\mathbbm{1}_{[0,100]}(d)$ together with
\begin{equation} 
\gamma_m = 
\begin{cases} 
10^{-1} \mathbbm{1}_{[0,2500]}(m) 
+ 10^{-2} \mathbbm{1}_{(2500,2750]}(m)
+ 10^{-3} \mathbbm{1}_{(2750,3000]}(m) & \mbox{ for } d\leq 100\\
10^{-1} \mathbbm{1}_{[0,1500]}(m) 
+ 10^{-2} \mathbbm{1}_{(1500,1750]}(m)
+ 10^{-3} \mathbbm{1}_{(1750,2000]}(m) & \mbox{ for }
d > 100.
\end{cases}
\end{equation} 
We set $\xi^{n,m,j}=(50,50,\dots,50) \in \R^{d}$ 
for every $(n,m,j) \in \{1,2, \dots, N\} \times \{0, 1, \dots, M-1\} \times \N$, and used 
feedforward neural networks with a $d$-dimensional 
input layer, two hidden layers of dimension $d+10+40\,\mathbbm{1}_{[1,100]}(d)$ and a one-dimensional output layer.

The reference values for $u(T,50,50,\dots,50)$ were calculated with the deep BSDE method 
of E et al.\ \cite{EHanJentzen2017}).
% and {\sc Matlab} code~ %\ref{code:branchingMatlab}
%in Appendix~ %\ref{subsec:BranchingMatlab} 
%below). 

%%%%%%%%%%%%%%%%%%%%%%%%%%%%%%%%%%%%%%%%%%%%%%%%%%%%%%%%%%%%%%%%%%%%%%%%%%
%%%%%%%%%%%%%%%%%%%%%%%%%%%%%%%%%%%%%%%%%%%%%%%%%%%%%%%%%%%%%%%%%%%%%%%%%
\subsection{Allen--Cahn-type equations}
\label{subsec:Allen-Cahn}
%%%
%FIGURE ALLEN-CAHN
%
%\begin{figure}
%	\includegraphics{AllenCahn50.pdf}
	%%%%%%%%%%%%%%%%%%%%%%%%%%%%%%%%%%
%	\caption{Plots of approximative 
%		calculations of the relative $ L^1 $-approximation error 
%		$\E\big[\tfrac{|\U^{\Theta_m}(\xi)-0.09909|}{0.09909}\big]$ 
%		and of the mean of the empirical loss function $\E\big[\tfrac{1}{J_m}
%		\sum_{j=1}^{J_m}|\Y^{m,\Theta_m,j,\bS_{m+1}}_N-g(\X^{m,j}_N)|^2\big]$ 
%		in the case of the $100$-dimensional 
%		Allen--Cahn equation \eqref{eq:example_pde_allen_cahn}
%		against $ m\in\{0,1,\ldots,2000\} $ .} 
%	\label{fig:figure_AllenCahn50}
	%%%%%%%%%%%%%%%%%%%%%%%%%%%%%%%%%%
%\end{figure}
%%%%%%%%%%%%%%%%%%%%%%%%%%%%%%%%%%%%%%%%%%5
Next, we approximate solutions of high-dimensional Allen--Cahn-type equations with a cubic nonlinearity of the form 
\begin{equation} \label{eq:AC}
 \tfrac{\partial}{\partial t}u(t,x) 
 = \Delta_x u(t,x)  + u(t,x) - [u(t,x)]^3, \quad(t,x) \in (0,T] \times \R^d,
\end{equation} 
with initial condition $u(0,x)=\arctan(\max_{i \in \{1,2,\dots,d\}} x_i)$; see also 
Beck et al.\ \cite[Section~4.1]{BeckEJentzen17}, E et al.\ \cite[Section~4.2]{EHanJentzen2017}, E et al.\ \cite[Section~3.4]{MultilevelPicard}, and Han et al.\ \cite[Section~3.3]{HanEJentzen17} for further numerical results for equation \eqref{eq:AC}. 

Table \ref{table:AC} lists deep splitting approximations of $u(T,0, \dots, 0)$ for $T = 0.3$, $N=10$, and different values of $d$.
We chose $\mu(x)=(0,0,\dots,0) \in \R^d$, $\sigma(x)=\sqrt{2}\mbox{ Id}_{\R^{d\times d}}$ and 
$f(x,y,z)= y-y^3$, $(x,y,z) \in \R^d \times \R \times \R^d$. We used $M= 500$ and
\begin{equation} 
\gamma_m = 10^{-1} \mathbbm{1}_{[0,300]}(m) 
+ 10^{-2} \mathbbm{1}_{(300,400]}(m)
+ 10^{-3} \mathbbm{1}_{(400,500]}(m).
\end{equation} 
We set $\xi^{n,m,j}=(0,0,\dots,0) \in \R^{d}$ for all $(n,m,j) \in \{1,2, \dots, N\} \times \{0, 1, \dots, M-1\} \times \N$
and used feedforward neural networks with a $d$-dimensional 
input layer, two hidden layers of dimension $d+10$ and a one-dimensional output layer.

%%%%%%%%%%%%%%%%%%%%%%%%%%%%%%%%%%%%%%%%%%%%%%
%table_AllenCahn50
%\input{table_AllenCahn50.tex}
%%%%%%%%%%%%%%%%%%%%%%%%
%%%%%%%%%%%%%%%%%%%%%%%%%%%%%%%%%%%%%%%%%%%%%%%%%%%%%%%%%%%%%
\begin{table}
	\begin{center}
		\resizebox{\textwidth}{!}{\begin{tabular}{|c|c|c|c|c|c|c|c|c|}
			\hline
			$ d $ & Expectation & Std.\ dev.\ & Ref.\ value & rel.\ $L^1$-error & Std.\ dev.\ rel.\ error & avg.\ runtime\\
			\hline
			10 & 0.89327 & 0.00299962 & 0.89060 & 0.00364 & 0.00258004 & 22.7s \\
			\hline
			50 & 1.01855 & 0.00073173 & 1.01830 & 0.00063 & 0.00036976 & 22.2s \\
			\hline
			100 & 1.04348 & 0.00029431 & 1.04510 & 0.00156 & 0.00028161 & 22.3s \\
			\hline
			200 & 1.06119 & 0.00018821 & 1.06220 & 0.00096 & 0.00017719 & 22.4s \\
			\hline
			300 & 1.06961 & 0.00017250 & 1.07217 & 0.00239 & 0.00016089 & 22.6s \\
			\hline
			500 & 1.07847 & 0.00013055 & 1.08124 & 0.00256 & 0.00012074 & 23.1s \\
			\hline
			1000 & 1.08842 & 0.00005689 & 1.09100 & 0.00236 & 0.00005215 & 25.9s \\
			\hline
			5,000 & 1.10522 & 0.00005201 & 1.10691 & 0.00153 & 0.00004699 & 134.6s \\
			\hline
			10,000 & 1.11071 & 0.00004502 & 1.11402 & 0.00296 & 0.00004041 & 473.6s \\
			\hline
		\end{tabular}}
			\caption{Deep splitting approximations of solutions of 
			Allen--Cahn-type equations of the form
			\eqref{eq:AC} for $T = 0.3$, $N = 10$ and different $d$.}
		\label{table:AC}
	\end{center}
\end{table}
%%%%%%%%%%%%%%%%%%%%%
%\noindent In Table~\ref{table:Allen-Cahn} %\ref{tab:table_AllenCahnPlainSGD.tex} 
%%%%%%%%%%%%%%%%%%%%%%%%%%%%%%%%%%%%
%we use 
%{\sc Python} code~\ref{code:AllenCahn} %%\ref{code:deepPDEmethodPlainSGDNoBN}
%in Subsection~\ref{subsec:Allen-Cahn} %%\ref{subsec:plainSGDCode} 
%below
%to 
%%%%%%%%%%%%%%%%%%%%%%%%%%%%%%%

The reference values for $u(T,0,0,\dots,0)$ were calculated with the multilevel Picard method;
%Branching diffusion method 
see, e.g.,  \cite{MLP2020, LinearScaling, MultilevelPicard,hutzenthaler2018overcoming,hutzenthaler2019overcoming,HutzenthalerKruse17}.
%%%%%%%%%%%%%%%%
% E et al.\ \cite{LinearScaling} and, e.g., Hutzenthaler et al.\ \cite{hutzenthaler2018overcoming}).
 %%%%%%%%%%%%%%%%%%%%%%%%%%%%%%%%%%%%%%%
% and {\sc Matlab} code~ %\ref{code:branchingMatlab}
%in Appendix~ %\ref{subsec:BranchingMatlab} 
%below).
%\\
%%%%%%%%%%%%%%%%%%%%%%%%%%%%%%%%%%%%%%%%%%%%%%%%%%%%%%%%%%%%%%%%%%
% 
% 
% 
%\\
\subsection{Semilinear heat equations}
\label{subsec:semiheat}

In this subsection, we calculate approximate solutions of semilinear heat equations of the form 
\begin{equation} \label{eq:HE}
\tfrac{\partial}{\partial t} u(t,x) = \Delta_x u(t,x) + \frac{1- |u(t,x)|^2}{1+ |u(t,x)|^2}, 
\quad (t,x) \in (0,T] \times \R^d,
\end{equation} 
with initial condition $u(0,x) = \nicefrac{5}{\left(10+{2}\|x\|_{\R^d}^2\right)}$.

Table \ref{table:HE} shows deep splitting approximations of $u(T,0, \dots, 0)$ for $T = 0.3$, $N = 20$, and 
different values of $d$. We derived them with $\mu(x)=(0,0,\dots,0) \in \R^d$, $\sigma(x)=\sqrt{2}\mbox{ Id}_{\R^{d\times d}}$,
and $f(x,y,z)= -\|z\|^2_{\R^d}$ for $(x,y,z) \in \R^d \times \R \times \R^d$. We used $M= 500$ and
\begin{equation}
\gamma_m = 10^{-1} \mathbbm{1}_{[0,300]}(m) 
+ 10^{-2} \mathbbm{1}_{(300,400]}(m)
+ 10^{-3} \mathbbm{1}_{(400,500]}(m).
\end{equation}
We set $\xi^{n,m,j}=(0,0,\dots,0) \in \R^{d}$ 
for all $(n,m,j) \in \{1,2, \dots, N\} \times \{0, 1, \dots, M-1\} \times \N$, and used feedforward neural networks with 
a $d$-dimensional input layer, two hidden layers of dimension $d+10$, and a one-dimensional output layer.
The reference values were computed with the multilevel Picard method; see, e.g.,  \cite{MLP2020, LinearScaling, MultilevelPicard,hutzenthaler2018overcoming,hutzenthaler2019overcoming,HutzenthalerKruse17}).  
%%%%%%%%%%%%%%%%%%%%%

To illustrate how the accuracy of the deep splitting method depends on the number of time steps $N$, 
we report in Table~\ref{table:HE-timesteps}, deep splitting approximations of $u(T,0, \dots, 0)$ for $d = 100$, $T = 0.3$, 
and different values of $N$. It can be seen that the decrease of the relative $L^1$-error is approximately linear in $N$.
For theoretical convergence results for the deep splitting method, we refer to Beck et al.\
\cite[Theorem~2]{Overview20}.
	
%Table~\ref{table:Nonlin-hiddenlayer}\!
%%%%%%%%%%%%%%%%%%%%%%%%%%%%%%%%%%
%%%%%%%%%%%%%%%%%%%%%%%%%%%%%%%%%
\begin{table}
	\begin{center}
		\resizebox{\textwidth}{!}{\begin{tabular}{|c|c|c|c|c|c|c|c|c|}
			\hline
			$ d $ & Expectation & Std.\ dev.\ & Ref.\ value & rel.\ $L^1$-error & Std.\ dev.\ rel.\ error & avg.\ runtime\\
			\hline
			10 & 0.47138 & 0.00035606 & 0.47006 & 0.00282 & 0.00075749 & 46.7s \\
			\hline
			50 & 0.34584 & 0.00018791 & 0.34425 & 0.00462 & 0.00054586 & 46.7s \\
			\hline
			100 & 0.31783 & 0.00008298 & 0.31674 & 0.00343 & 0.00026198 & 47.4s \\
			\hline
			200 & 0.30210 & 0.00002238 & 0.30091 & 0.00394 & 0.00007436 & 48.1s \\
			\hline
			300 & 0.29654 & 0.00001499 & 0.29534 & 0.00406 & 0.00005075 & 48.3s \\
			\hline
			500 & 0.29200 & 0.00000611 & 0.29095 & 0.00361 & 0.00002099 & 48.5s \\
			\hline
			1000 & 0.28852 & 0.00000267 & 0.28753 & 0.00344 & 0.00000930 & 54.1s \\
			\hline
			5,000 & 0.28569 & 0.00000042 & 0.28469 & 0.00352 & 0.00000148 & 286.3s \\
			\hline
			10,000 & 0.28533 & 0.00000048 & 0.28433 & 0.00353 & 0.00000170 & 1013.0s \\
			\hline
		\end{tabular}}
	\caption{Deep splitting approximations of solutions of the semilinear heat equation
	\eqref{eq:HE} for $T = 0.3$, $N=20$ and different $d$.}
	\label{table:HE}
	\end{center}
\end{table}
%%%%%%%%%%%%%%%%%%%%%%%%%%%%%%%%%%%%%%%%%%%%%%
%
%%%%%%%%%%%%%%%%%%%%%%%%%%%%%%%%%
\begin{table}
	\begin{center}
		\resizebox{\textwidth}{!}{\begin{tabular}{|c|c|c|c|c|c|c|c|c|}
				\hline
				$ N $ & Expectation & Std.\ dev.\ & Ref.\ value & rel.\ $L^1$-error & Std.\ dev.\ rel.\ error & avg.\ runtime\\
				\hline
				1 & 0.33820 & 0.00005021 & 0.31674 & 0.06777 & 0.00015852 & 0.5s \\
				\hline
				2 & 0.32836 & 0.00006189 & 0.31674 & 0.03669 & 0.00019539 & 2.2s \\
				\hline
				4 & 0.32259 & 0.00003834 & 0.31674 & 0.01848 & 0.00012104 & 5.7s \\
				\hline
				8 & 0.31969 & 0.00006929 & 0.31674 & 0.00931 & 0.00021878 & 12.8s \\
				\hline
				16 & 0.31810 & 0.00006968 & 0.31674 & 0.00430 & 0.00021998 & 27.1s \\
				\hline
				20 & 0.31783 & 0.00010565 & 0.31674 & 0.00345 & 0.00033354 & 34.3s \\
				\hline
				32 & 0.31739 & 0.00006564 & 0.31674 & 0.00206 & 0.00020723 & 55.6s \\
				\hline
		\end{tabular}}
		\caption{Deep splitting approximations of solutions of the semilinear heat equation
	\eqref{eq:HE} for $d=100$, $T = 0.3$ and different $N$.}
		\label{table:HE-timesteps}
	\end{center}
\end{table}
\subsection{Sine-Gordon-type equations}
\label{subsec:Sine-Gordon}
As a last example, we approximate solutions of sine-Gordon-type equations of the form 
\begin{equation} \label{eq:SG}
\tfrac{\partial}{\partial t} u(t,x) = \Delta_x u(t,x) + \sin(u(t,x)), \quad 
(t,x) \in (0,T] \times \R^d, 
\end{equation} 
with initial condition $u(0,x) = \nicefrac{5}{\left(10+{2}\|x\|_{\R^d}^2 \right)}$.

Table \ref{table:SG} shows deep splitting approximations of $u(T, 0, \dots, 0)$ for $T = 0.3$, $N = 20$
and different values of $d$.
We chose $\mu(x)=(0,0,\dots,0) \in \R^d$, $\sigma(x)=\sqrt{2}\mbox{ Id}_{\R^{d\times d}}$, and
$f(x,y,z)= \sin(y)$, $(x,y,z) \in \R^d \times \R \times \R^d$. 
We used $M= 1000$ and
\begin{equation}
\gamma_{m} = 10^{-1} \mathbbm{1}_{[0,250]}(m) 
+ 10^{-2} \mathbbm{1}_{(250,500]}(m)
+ 10^{-3} \mathbbm{1}_{(500,750]}(m)
+ 10^{-4} \mathbbm{1}_{(750,1000]}(m).
\end{equation} 
We set $\xi^{n,m,j}=(0,0,\dots,0) \in \R^{d}$ for every $(n,m,j) \in \{1,2, \dots, N\} \times \{0, 1, \dots, M-1\} \times \N$
and used feedforward neural networks with a $d$-dimensional 
input layer, two hidden layers of dimension $d+50$ and a one-dimensional output layer.

The reference values for $u(T,0,0,\dots,0)$ were computed with the multilevel Picard method;
%Branching diffusion method 
%%%%%%%%%%%%%%%%
%%%% FALLS MULTI-PICARD
see, e.g., \cite{ LinearScaling, MultilevelPicard,hutzenthaler2018overcoming,hutzenthaler2019overcoming,HutzenthalerKruse17}. 
%

%%%%%%%%%%%%%%%%%%%%%%%%%%%%%%%%%%%%%%%%%%%%%%%%%%%%%%%%%%%%%%%
\begin{table}
	\begin{center}
		\resizebox{\textwidth}{!}{\begin{tabular}{|c|c|c|c|c|c|c|c|c|}
				\hline
				$ d $ & Expectation & Std.\ dev.\ & Ref.\ value & rel.\ $L^1$-error & Std.\ dev.\ rel.\ error & avg.\ runtime\\
				\hline
				10 & 0.3218822 & 0.00069331 & 0.3229470 & 0.0032972 & 0.0021468 & 55.6s \\
				\hline
				50 & 0.0990598 & 0.00013433 & 0.0993633 & 0.0030544 & 0.0013519 & 55.4s \\
				\hline
				100 & 0.0526955 & 0.00005390 & 0.0528368 & 0.0026741 & 0.0010202 & 55.6s \\
				\hline
				200 & 0.0271860 & 0.00001617 & 0.0272410 & 0.0020176 & 0.0005936 & 56.2s \\
				\hline
				300 & 0.0183162 & 0.00001010 & 0.0183617 & 0.0024765 & 0.0005502 & 55.8s \\
				\hline
				500 & 0.0110819 & 0.00000428 & 0.0111071 & 0.0022647 & 0.0003851 & 54.7s \\
				\hline
				1000 & 0.0055775 & 0.00000090 & 0.0055896 & 0.0021697 & 0.0001613 & 58.1s \\
				\hline
				5,000 & 0.0011209 & 0.00000012 & 0.0011231 & 0.0019464 & 0.0001070 & 332.5s \\
				\hline
				10,000 & 0.0005608 & 0.00000004 & 0.0005621 & 0.0022495 & 0.0000720 & 1083.4s \\
				\hline
		\end{tabular}}
		\caption{Deep splitting approximations of the solution of the 
			sine-Gordon-type equation \eqref{eq:SG} for $T= 0.3$, $N=20$ and different $d$.}
			\label{table:SG}
	\end{center}
\end{table}
%%%%%%%%%%%%%%%%%%%%%%%%%%%%%%%%%%%%%%%%%%%%%%%%%%%%%%%%%%%
%
%
\section{Comparison with other methods}
\label{subsec:comparison}

In this section, we compare the deep splitting method with the deep BSDE method of 
E et al.\ \cite{EHanJentzen2017} and the multilevel Picard method of \cite{MLP2020,hutzenthaler2018overcoming},
which also have been shown to produce good results in approximating solutions of high-dimensional PDEs.

\begin{figure}[ht!]
\centering
	{%
\includegraphics[width=.75\textwidth]{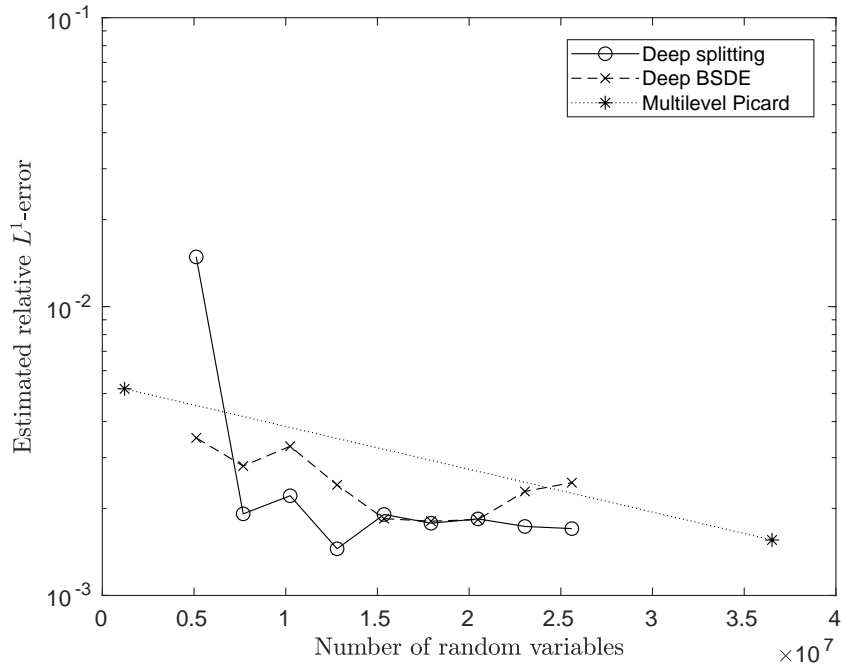} 
		% 0.92 %0.65
	}
\caption{Estimated relative $L^1$-errors as a function of the number of one-dimensional standard normal 
random variables used by the deep splitting algorithm, the deep BSDE approach of E et al.\ \cite{EHanJentzen2017} and 
the multilevel Picard method of \cite{MLP2020,hutzenthaler2018overcoming} 
for the sine-Gordon-type equation \eqref{eq:SG} with $d=10$ and $T = 0.3$.}
\label{fig:Comparison-Deep}
\end{figure}

Figure~\ref{fig:Comparison-Deep} shows estimated relative $L^1$-approximation errors as a function
of the number of one-dimensional standard normal random variables used by the three methods to approximate 
$u(T, 0, \dots, 0)$ for the solution $u \colon [0,T] \times \R^d \to \R$ of the sine-Gordon-type equation \eqref{eq:SG} 
for $d = 10$ and $T = 0.3$. It can be seen that for this particular example, the three methods yield 
comparable results. But it has to be 
noted that the deep splitting algorithm and the deep BSDE method both involve different 
hyper-parameters, which, for good results, have to be fine-tuned depending on the 
form and the parameters of the PDE. Figure~\ref{fig:Comparison-Deep} just shows 
results for particular implementations of the three methods. 

Generally, since the deep splitting method uses neural networks to approximate the solution 
$u \colon [0,T] \times \R^d \to \R$ of a PDE on a time grid $0 = t_0 < t_1 < \dots < t_N = T$, 
it can learn approximations of $u(t_n, x)$, $n \in \{1, 2, \dots, N\}$, simultaneously for all $x \in \R^d$.
Similarly, the deep BSDE method can be implemented so that it approximates $u(T,x)$ directly for all $x \in \R^d$.
However, for a temporal discretization with $N$ subintervals, it then needs to train $N$ neural networks 
at the same time, whereas the deep splitting method trains one network after the other. 
So even if the two approaches need a similar total number of one-dimensional standard normal random variables
to achieve a given accuracy, the deep splitting method can handle larger problems since it 
decomposes them into smaller computational tasks which can be solved successively.

The advantage of the multilevel Picard method is that there exist theoretical bounds on the
computational effort needed for a given approximation accuracy; 
see e.g., \cite{LinearScaling, hutzenthaler2018overcoming, hutzenthaler2019overcoming, HutzenthalerKruse17}. 
However, the method needs to calculate approximations of $u(t,x)$ for different space-time points
$(t, x)$ separately and becomes impractical for large $t$.

%A more comprehensive comparison between the deep splitting approximation method, the deep BSDE approximation method, and the MLP approximation method for high-dimensional parabolic PDEs is left %remains a  topic
%for future research.
%Therefore, it is possible that the deep splitting approximation method may outperform the multilevel Picard approximation method when considering different PDEs.
%In addition, note that the advantage of the deep splitting approximation method, compared to the multilevel Picard approximation method, is that it can be applied to PDEs which depend non-linearly on the gradient, as well
%as to a general class of coefficients in the PDE. Furthermore, one may extend the deep splitting approximation method to, e.g., free boundary problems.

\vspace{-0.1cm}
%%%%%%%%%%%%%%%%%%%%%%%%%%%%%%%%%%%%%%%%%%%%%%%%%%%%%%%%%%%%

\section{Conclusion}

In this paper we have developed a new numerical method to approximate solutions 
of high-dimensional nonlinear parabolic PDEs. It splits the differential operator into a 
linear and a nonlinear part and uses deep learning together with the 
Feynman--Kac formula to iteratively solve linear approximations of the equation over small time intervals.
This breaks the PDE approximation task into smaller problems that can be solved successively.
As a consequence, the approach can be applied to extremely high-dimensional nonlinear 
PDEs. We have tested the method on Hamilton--Jacobi--Bellman equations,
nonlinear Black--Scholes equations, Allen--Cahn-type equations, semilinear heat equations as well as
sine-Gordon-type equations. In all cases, it has produced accurate results in high dimensions
with short run times.

%%%%%%%%%%%%%%%%%%%%%
%%%%%%%%%%%%%%%%%%%%%%%%%%%%%%%%%%%%%%%%%%%%%%%%%%%%%%%%%%%%%%%%%%%%%%%%%%%
%
%\newpage
%%%%%%%%%%%%%%%%%%%%%%%%%%%%%%%%%%%%%%%%%%%%%%%%%%%%%%%%%%%%%
 \subsection*{Acknowledgments}
We are grateful to Adam Andersson for fruitful discussions.
This project has been supported by the Swiss National Science Foundation
Grant $200020\_175699$  ``Higher order numerical approximation methods 
for stochastic partial differential equations'', by the Deutsche Forschungsgemeinschaft 
under Germany’s Excellence Strategy EXC 2044-390685587, Mathematics M\"unster: 
Dynamics - Geometry - Structure and by the Nanyang Assistant Professorship Grant
``Machine Learning based Algorithms in Finance and Insurance''.
%%%%%%%%%%%%%%%%%%%%%%%%%%%%%%%
\bibliographystyle{acm}
%\bibliography{bibfile}

\begin{thebibliography}{10}
	
	\bibitem{BallyPages2003}
	{\sc Bally, V., and Pag\`es, G.}
	\newblock A quantization algorithm for solving multi-dimensional discrete-time
	optimal stopping problems.
	\newblock {\em Bernoulli 9}, 6 (2003), 1003--1049.
	
	\bibitem{DeepKolmogorov}
	{\sc Beck, C., Becker, S., Grohs, P., Jaafari, N., and Jentzen, A.}
	\newblock Solving stochastic differential equations and {K}olmogorov equations
	by means of deep learning.
	\newblock {\em Revision requested from Journal of Scientific Computing.
		arXiv:1806.00421\/} (2018), 56~pages.
	
	\bibitem{BeckEJentzen17}
	{\sc Beck, C., E, W., and Jentzen, A.}
	\newblock Machine learning approximation algorithms for high-dimensional fully
	nonlinear partial differential equations and second-order backward stochastic
	differential equations.
	\newblock {\em Journal of Nonlinear Science\/} (2017), 1--57.
	

\bibitem{Overview20}
{\sc Beck, C., Hutzenthaler, M., Jentzen, A., and Kuckuck, B.}
\newblock An overview on deep learning-based approximation
methods for partial differential equations.
\newblock {\em arXiv:2012.12348\/} (2020), 22~pages.

%\newpage
	
\bibitem{MLP2020}
	{\sc Becker, S., Braunwarth, R., Hutzenthaler, M., Jentzen, A., and von Wurstemberger, P.}
	\newblock Numerical simulations for full history recursive multilevel {P}icard approximations for systems of high-dimensional partial differential equations.
	\newblock {\em Accepted in Communications in Computational Physics. arXiv:2005.10206\/} (2020), 21~pages.
	
\bibitem{becker2018deep}
	{\sc Becker, S., Cheridito, P., and Jentzen, A.}
	\newblock Deep optimal stopping.
	\newblock {\em The {J}ournal of {M}achine {L}earning {R}esearch 20}, 74 (2015),
	1--25.
	
	\bibitem{BenderDenk2007}
	{\sc Bender, C., and Denk, R.}
	\newblock A forward scheme for backward {SDE}s.
	\newblock {\em Stochastic Processes and their Applications 117}, 12 (2007),
	1793--1812.
	
	\bibitem{Bender2015Primal}
	{\sc Bender, C., Schweizer, N., and Zhuo, J.}
	\newblock A primal-dual algorithm for {BSDEs}.
	\newblock {\em Mathematical Finance 27}, 3 (2017), 866--901.
	
	\bibitem{Bengio09}
	{\sc Bengio, Y.}
	\newblock Learning deep architectures for {AI}.
	\newblock {\em Foundations and Trends in Machine Learning 2}, 1 (2009), 1--127.
	
	\bibitem{berg2018unified}
	{\sc Berg, J., and Nystr{\"o}m, K.}
	\newblock A unified deep artificial neural network approach to partial
	differential equations in complex geometries.
	\newblock {\em Neurocomputing 317\/} (2018), 28--41.
	
	\bibitem{berner2018analysis}
	{\sc Berner, J., Grohs, P., and Jentzen, A.}
	\newblock Analysis of the generalization error: Empirical risk minimization
	over deep artificial neural networks overcomes the curse of dimensionality in
	the numerical approximation of {B}lack-{S}choles partial differential
	equations.
	\newblock {\em SIAM Journal on Mathematics of Data Science 2}, 3 (2020), 631--657.
	
%		\bibitem{berner2018analysis}
%	{\sc Berner, J., Grohs, P., and Jentzen, A.}
%	\newblock Analysis of the generalization error: Empirical risk minimization
%	over deep artificial neural networks overcomes the curse of dimensionality in
%	the numerical approximation of {B}lack-{S}choles partial differential
%	equations.
%	\newblock {\em arXiv:1809.03062\/} (2018), 35 pages.
	
	\bibitem{BouchardTouzi2004}
	{\sc Bouchard, B., and Touzi, N.}
	\newblock Discrete-time approximation and {M}onte-{C}arlo simulation of
	backward stochastic differential equations.
	\newblock {\em Stochastic Processes and their Applications 111}, 2 (2004),
	175--206.
	
	\bibitem{Braess2007FEM}
	{\sc Braess, D.}
	\newblock {\em Finite elements}, third~ed.
	\newblock Cambridge University Press, Cambridge, 2007.
	\newblock Theory, fast solvers, and applications in elasticity theory,
	Translated from the German by Larry L. Schumaker.
	
	\bibitem{chan2018machine}
	{\sc Chan-Wai-Nam, Q., Mikael, J., and Warin, X.}
	\newblock Machine learning for semi~linear {PDEs}.
	\newblock {\em Journal of Scientific Computing 79}, 3 (2019), 1667--1712.
	
%		\bibitem{chan2018machine}
%	{\sc Chan-Wai-Nam, Q., Mikael, J., and Warin, X.}
%	\newblock Machine learning for semi~linear {PDEs}.
%	\newblock {\em arXiv:1809.07609\/} (2018), 38 pages.
	
	\bibitem{Chassagneux2014}
	{\sc Chassagneux, J.-F.}
	\newblock Linear multistep schemes for {BSDE}s.
	\newblock {\em SIAM Journal on Numerical Analysis 52}, 6 (2014), 2815--2836.
	
	\bibitem{ChassagneuxCrisan2014}
	{\sc Chassagneux, J.-F., and Crisan, D.}
	\newblock Runge-{K}utta schemes for backward stochastic differential equations.
	\newblock {\em The Annals of Applied Probability 24}, 2 (2014), 679--720.
	
	\bibitem{ChassagneuxRichou2015}
	{\sc Chassagneux, J.-F., and Richou, A.}
	\newblock Numerical stability analysis of the {E}uler scheme for {BSDE}s.
	\newblock {\em SIAM Journal on Numerical Analysis 53}, 2 (2015), 1172--1193.
	
	\bibitem{ChassagneuxRichou2016}
	{\sc Chassagneux, J.-F., and Richou, A.}
	\newblock Numerical simulation of quadratic {BSDE}s.
	\newblock {\em The Annals of Applied Probability 26}, 1 (2016), 262--304.
	
	\bibitem{CrisanManolarakis2010}
	{\sc Crisan, D., and Manolarakis, K.}
	\newblock Probabilistic methods for semilinear partial differential equations.
	{A}pplications to finance.
	\newblock {\em M2AN Mathematical Modelling and Numerical Analysis 44}, 5
	(2010), 1107--1133.
	
	\bibitem{CrisanManolarakis2012}
	{\sc Crisan, D., and Manolarakis, K.}
	\newblock Solving backward stochastic differential equations using the cubature
	method: application to nonlinear pricing.
	\newblock {\em SIAM Journal on Financial Mathematics 3}, 1 (2012), 534--571.
	
	\bibitem{CrisanManolarakis2014}
	{\sc Crisan, D., and Manolarakis, K.}
	\newblock Second order discretization of backward {SDE}s and simulation with
	the cubature method.
	\newblock {\em The Annals of Applied Probability 24}, 2 (2014), 652--678.
	
	\bibitem{CrisanManolarakisTouzi2010}
	{\sc Crisan, D., Manolarakis, K., and Touzi, N.}
	\newblock On the {M}onte {C}arlo simulation of {BSDE}s: an improvement on the
	{M}alliavin weights.
	\newblock {\em Stochastic Processes and their Applications 120}, 7 (2010),
	1133--1158.
	
	\bibitem{DeckKruse_ParametrixMethod2002}
	{\sc Deck, T., and Kruse, S.}
	\newblock Parabolic differential equations with unbounded coefficients---a
	generalization of the parametrix method.
	\newblock {\em Acta Applicandae Mathematicae 74}, 1 (2002), 71--91.
	
	\bibitem{DelarueMenozzi2006}
	{\sc Delarue, F., and Menozzi, S.}
	\newblock A forward-backward stochastic algorithm for quasi-linear {PDE}s.
	\newblock {\em The Annals of Applied Probability 16}, 1 (2006), 140--184.
	
	\bibitem{Dorsek12}
	{\sc D{\"o}rsek, P.}
	\newblock Semigroup splitting and cubature approximations for the stochastic
	{N}avier-{S}tokes equations.
	\newblock {\em SIAM Journal on Numerical Analysis 50}, 2 (2012), 729--746.
	
	\bibitem{DouglasMaProtter}
	{\sc Douglas, Jr., J., Ma, J., and Protter, P.}
	\newblock Numerical methods for forward-backward stochastic differential
	equations.
	\newblock {\em The Annals of Applied Probability 6}, 3 (1996), 940--968.
	
	\bibitem{DuffieEtAl1996DefaultableSecurities}
	{\sc Duffie, D., Schroder, M., and Skiadas, C.}
	\newblock Recursive valuation of defaultable securities and the timing of
	resolution of uncertainty.
	\newblock {\em Ann. Appl. Probab. 6}, 4 (1996), 1075--1090.
	
	\bibitem{EHanJentzen2017}
	{\sc {E}, W., {Han}, J., and {Jentzen}, A.}
	\newblock Deep learning-based numerical methods for high-dimensional parabolic
	partial differential equations and backward stochastic differential equations.
	\newblock {\em Communications in Mathematics and Statistics 5\/} (2017),
	349--380.
	
	\bibitem{LinearScaling}
	{\sc E, W., Hutzenthaler, M., Jentzen, A., and Kruse, T.}
	\newblock Multilevel {P}icard iterations for solving smooth semilinear
	parabolic heat equations.
	\newblock {\em Revision requested from SN Partial Differential Equations and Applications. arXiv:1607.03295\/} (2017), 18~pages.
	
%	\bibitem{LinearScaling}
%	{\sc E, W., Hutzenthaler, M., Jentzen, A., and Kruse, T.}
%	\newblock Multilevel {P}icard iterations for solving smooth semilinear
%	parabolic heat equations.
%	\newblock {\em arXiv:1607.03295\/} (2017), 18 pages.
	
	\bibitem{MultilevelPicard}
	{\sc E, W., Hutzenthaler, M., Jentzen, A., and Kruse, T.}
	\newblock On multilevel {P}icard numerical approximations for high-dimensional
	nonlinear parabolic partial differential equations and high-dimensional
	nonlinear backward stochastic differential equations.
	\newblock {\em Journal of Scientific Computing\/} (2019), 1--38.
	
	\bibitem{EYu17}
	{\sc E, W., and Yu, B.}
	\newblock The deep {R}itz method: A deep learning-based numerical algorithm for
	solving variational problems.
	\newblock {\em Communications in Mathematics and Statistics 6}, 1 (2018), 1--12.
	
%	\bibitem{EYu17}
%	{\sc E, W., and Yu, B.}
%	\newblock The deep {R}itz method: A deep learning-based numerical algorithm for
%	solving variational problems.
%	\newblock {\em arXiv:1710.00211\/} (2017), 14 pages.
	
	\bibitem{elbrachter2018dnn}
	{\sc Elbr{\"a}chter, D., Grohs, P., Jentzen, A., and Schwab, C.}
	\newblock {DNN} expression rate analysis of high-dimensional {PDEs}:
	Application to option pricing.
	\newblock {\em 
		Accepted in Constructive Approximation. arXiv:1809.07669\/} (2018), 50 pages.
	
	\bibitem{FarahmandNabiNikovski17}
	{\sc Farahmand, A.-m., Nabi, S., and Nikovski, D.~N.}
	\newblock Deep reinforcement learning for partial differential equation
	control.
	\newblock {\em 2017 American Control Conference (ACC)\/} (2017), 3120--3127.
	
	\bibitem{FujiiTakahashiTakahashi17}
	{\sc Fujii, M., Takahashi, A., and Takahashi, M.}
	\newblock Asymptotic expansion as prior knowledge in deep learning method for
	high dimensional {BSDEs}.
	\newblock {\em Asia-Pacific Financial Markets 26}, 3 (2019), 391--408.
	
%		\bibitem{FujiiTakahashiTakahashi17}
%	{\sc Fujii, M., Takahashi, A., and Takahashi, M.}
%	\newblock Asymptotic expansion as prior knowledge in deep learning method for
%	high dimensional {BSDEs}.
%	\newblock {\em arXiv:1710.07030\/} (2017), 16 pages.
	
	\bibitem{glorot2010understanding}
	{\sc Glorot, X., and Bengio, Y.}
	\newblock Understanding the difficulty of training deep feedforward neural
	networks.
	\newblock {\em Proceedings of the thirteenth international conference on
		artificial intelligence and statistics\/} (2010), 249--256.
	
	\bibitem{GobetLabart2010}
	{\sc Gobet, E., and Labart, C.}
	\newblock Solving {BSDE} with adaptive control variate.
	\newblock {\em SIAM Journal on Numerical Analysis 48}, 1 (2010), 257--277.
	
	\bibitem{GobetLemor2008Numerical}
	{\sc Gobet, E., and Lemor, J.-P.}
	\newblock Numerical simulation of {BSDE}s using empirical regression methods:
	theory and practice.
	\newblock {\em arXiv:0806.4447\/} (2008), 17 pages.
	
	\bibitem{GobetLemorWarin2005}
	{\sc Gobet, E., Lemor, J.-P., and Warin, X.}
	\newblock A regression-based {M}onte {C}arlo method to solve backward
	stochastic differential equations.
	\newblock {\em The Annals of Applied Probability 15}, 3 (2005), 2172--2202.
	
	\bibitem{GobetLopezSalasTurkedjiev2016}
	{\sc Gobet, E., L\'opez-Salas, J.~G., Turkedjiev, P., and V\'azquez, C.}
	\newblock Stratified regression {M}onte-{C}arlo scheme for semilinear {PDE}s
	and {BSDE}s with large scale parallelization on {GPU}s.
	\newblock {\em SIAM Journal on Scientific Computing 38}, 6 (2016), C652--C677.
	
	\bibitem{GobetTurkedjiev2016}
	{\sc Gobet, E., and Turkedjiev, P.}
	\newblock Approximation of backward stochastic differential equations using
	{M}alliavin weights and least-squares regression.
	\newblock {\em Bernoulli 22}, 1 (2016), 530--562.
	
	\bibitem{GobetTurkedjiev2016MathComp}
	{\sc Gobet, E., and Turkedjiev, P.}
	\newblock Linear regression {MDP} scheme for discrete backward stochastic
	differential equations under general conditions.
	\newblock {\em Mathematics of Computation 85}, 299 (2016), 1359--1391.
	
	\bibitem{goudenege2019machine}
	{\sc Goudenege, L., Molent, A., and Zanette, A.}
	\newblock Machine learning for pricing {A}merican options in high-dimensional {M}arkovian and non-{M}arkovian models
	\newblock {\em Quantitative Finance 20\/} 4 (2020), 573--591.
	
	\bibitem{GrekschLisei_ApproximationOfStochasticNonlinearEquationsBySplittingMethod2013}
	{\sc Grecksch, W., and Lisei, H.}
	\newblock Approximation of stochastic nonlinear equations of {S}chr\"odinger
	type by the splitting method.
	\newblock {\em Stochastic Analysis and Applications 31}, 2 (2013), 314--335.
	
	\bibitem{grohs2018proof}
	{\sc Grohs, P., Hornung, F., Jentzen, A., and von Wurstemberger, P.}
	\newblock A proof that artificial neural networks overcome the curse of
	dimensionality in the numerical approximation of {B}lack-{S}choles partial
	differential equations.
	\newblock {\em Accepted in the Memoirs of the American Mathematical Society. arXiv:1809.02362\/} (2018), 124 pages.
	
	\bibitem{GyoengyKrylov_OnTheRateOfConvergenceOfSplittingUpApproximationsForSPDEs2003}
	{\sc Gy{\"o}ngy, I., and Krylov, N.}
	\newblock On the rate of convergence of splitting-up approximations for
	{SPDEs}.
	\newblock {\em Stochastic inequalities and applications\/} (2003), 301--321.
	
	\bibitem{GyoengyKrylov_OnTheSplittingUpMethodForSPDEs}
	{\sc Gy\"ongy, I., and Krylov, N.}
	\newblock On the splitting-up method and stochastic partial differential
	equations.
	\newblock {\em The Annals of Probability 31}, 2 (2003), 564--591.
	
	\bibitem{HairerHutzenthalerJentzen_LossOfRegularity2015}
	{\sc Hairer, M., Hutzenthaler, M., and Jentzen, A.}
	\newblock Loss of regularity for {K}olmogorov equations.
	\newblock {\em The Annals of Probability 43}, 2 (2015), 468--527.
	
	\bibitem{HanEJentzen17}
	{\sc Han, J., Jentzen, A., and E, W.}
	\newblock Solving high-dimensional partial differential equations using deep
	learning.
	\newblock {\em Proceedings of the National Academy of Sciences 115}, 34 (2018),
	8505--8510.
	
	\bibitem{han2018convergence}
	{\sc Han, J., and Long, J.}
	\newblock Convergence of the deep {BSDE} method for coupled {FBSDEs}.
	\newblock {\em Probability, Uncertainty and Quantitative Risk 5\/} 1 (2020), 1--33.
	
	\bibitem{Labordere2012}
	{\sc Henry-Labord{\`e}re, P.}
	\newblock Counterparty risk valuation: a marked branching diffusion approach.
	\newblock {\em arXiv:1203.2369\/} (2012), 17 pages.
	
	\bibitem{HenryLabordere17}
	{\sc Henry-Labordere, P.}
	\newblock Deep primal-dual algorithm for {BSDEs}: Applications of machine
	learning to {CVA} and {IM}.
	\newblock {\em Preprint, SSRN--id3071506\/} (2017), 16 pages.
	
	\bibitem{Labordereetal2016arxiv}
	{\sc Henry-Labord{\`e}re, P., Oudjane, N., Tan, X., Touzi, N., and Warin, X.}
	\newblock Branching diffusion representation of semilinear {PDE}s and {M}onte
	{C}arlo approximation.
	\newblock {\em {A}nnales de {l'I}nstitut {H}enri {P}oincar{\'e}
		{(B)} {P}robabilités et {S}tatistiques 55}, 1 (2019), 184--210.
	
	
%		\bibitem{Labordereetal2016arxiv}
%	{\sc Henry-Labord{\`e}re, P., Oudjane, N., Tan, X., Touzi, N., and Warin, X.}
%	\newblock Branching diffusion representation of semilinear {PDE}s and {M}onte
%	{C}arlo approximation.
%	\newblock {\em To appear in {A}nnales de {l'I}nstitut {H}enri {P}oincar{\'e}
%		{(B)} {P}robabilités et {S}tatistiques, arXiv:1603.01727\/} (2016), 30
%	pages.
	
	\bibitem{LabordereTanTouzi2014}
	{\sc Henry-Labord{\`e}re, P., Tan, X., and Touzi, N.}
	\newblock A numerical algorithm for a class of {BSDE}s via the branching
	process.
	\newblock {\em Stochastic Processes and their Applications 124}, 2 (2014),
	1112--1140.
	
	\bibitem{HuijskensRuijterOosterlee2016}
	{\sc Huijskens, T.~P., Ruijter, M.~J., and Oosterlee, C.~W.}
	\newblock Efficient numerical {F}ourier methods for coupled forward-backward
	{SDE}s.
	\newblock {\em Journal of Computational and Applied Mathematics 296\/} (2016),
	593--612.
	
	\bibitem{hure2019some}
	{\sc Hur{\'e}, C., Pham, H., and Warin, X.}
	\newblock Some machine learning schemes for high-dimensional nonlinear {PDEs}.
	\newblock {\em arXiv:1902.01599\/} (2019), 33 pages.
	
	\bibitem{hutzenthaler2019proof}
	{\sc Hutzenthaler, M., Jentzen, A., Kruse, T., and Nguyen, T.~A.}
	\newblock A proof that rectified deep neural networks overcome the curse of
	dimensionality in the numerical approximation of semilinear heat equations.
	\newblock {\em SN Partial Differential Equations and Applications 1\/} (2020), 1--34.
	
	\bibitem{hutzenthaler2018overcoming}
	{\sc Hutzenthaler, M., Jentzen, A., Kruse, T., Nguyen, T.~A., and von
		Wurstemberger, P.}
	\newblock Overcoming the curse of dimensionality in the numerical approximation
	of semilinear parabolic partial differential equations.
	\newblock {\em Accepted in Proceedings of the Royal Society of London. Series A. arXiv:1807.01212\/} (2018), 27 pages.
	
	\bibitem{HutzenthalerJentzenSalimova16}
	{\sc Hutzenthaler, M., Jentzen, A., and Salimova, D.}
	\newblock Strong convergence of full-discrete nonlinearity-truncated
	accelerated exponential {E}uler-type approximations for stochastic
	{K}uramoto-{S}ivashinsky equations.
	\newblock {\em Communications in Mathematical Sciences 16\/} (2018),
	1489--1529.
	
	\bibitem{hutzenthaler2019overcoming}
	{\sc Hutzenthaler, M., Jentzen, A., and von Wurstemberger, P.}
	\newblock Overcoming the curse of dimensionality in the approximative pricing
	of financial derivatives with default risks.
	\newblock {\em Electronic Journal of Probability 25\/}, (2020), 1--73.
	
	\bibitem{HutzenthalerKruse17}
	{\sc Hutzenthaler, M., and Kruse, T.}
	\newblock Multi-level {P}icard approximations of high-dimensional semilinear
	parabolic differential equations with gradient-dependent nonlinearities.
	\newblock {\em SIAM Journal on Numerical Analysis 58}, 2 (2020), 929--961.
	
%		\bibitem{HutzenthalerKruse17}
%	{\sc Hutzenthaler, M., and Kruse, T.}
%	\newblock Multi-level {P}icard approximations of high-dimensional semilinear
%	parabolic differential equations with gradient-dependent nonlinearities.
%	\newblock {\em arXiv:1711.01080\/} (2017), 19 pages.
	
	\bibitem{IoffeSzegedy2015}
	{\sc Ioffe, S., and Szegedy, C.}
	\newblock Batch {N}ormalization: {A}ccelerating {D}eep {N}etwork {T}raining by
	{R}educing {I}nternal {C}ovariate {S}hift.
	\newblock {\em arXiv:1502.03167\/} (2015), 11 pages.
	
	\bibitem{jacquier2019deep}
	{\sc Jacquier, A., and Oumgari, M.}
	\newblock Deep {PPDEs} for rough local stochastic volatility.
	\newblock {\em arXiv:1906.02551\/} (2019), 21 pages.
	
	\bibitem{jentzen2018proof}
	{\sc Jentzen, A., Salimova, D., and Welti, T.}
	\newblock A proof that deep artificial neural networks overcome the curse of
	dimensionality in the numerical approximation of {K}olmogorov partial
	differential equations with constant diffusion and nonlinear drift
	coefficients.
	\newblock {\em Accepted in Communications in Mathematical Sciences. arXiv:1809.07321\/} (2018), 48 pages.
	
	\bibitem{KingmaBa2015}
	{\sc Kingma, D., and Ba, J.}
	\newblock {Adam: {A} {M}ethod for {S}tochastic {O}ptimization}.
	\newblock {\em arXiv:1412.6980\/} (2014), 15 pages.
	
	\bibitem{Klenke_2014}
	{\sc Klenke, A.}
	\newblock {\em Probability theory. A comprehensive course}, second~ed.
	\newblock Universitext. Springer, London, 2014.
	
	\bibitem{KloedenPlaten1992}
	{\sc Kloeden, P.~E., and Platen, E.}
	\newblock {\em Numerical solution of stochastic differential equations},
	vol.~23 of {\em Applications of Mathematics (New York)}.
	\newblock Springer-Verlag, Berlin, 1992.
	
	\bibitem{Krylov_LecturesHoelder1996}
	{\sc Krylov, N.~V.}
	\newblock {\em Lectures on elliptic and parabolic equations in {H}\"older
		spaces}, vol.~12 of {\em Graduate Studies in Mathematics}.
	\newblock American Mathematical Society, Providence, RI, 1996.
	
	\bibitem{Krylov_KolmogorovsEquations1998}
	{\sc Krylov, N.~V.}
	\newblock On {K}olmogorov's equations for finite-dimensional diffusions.
	\newblock {\em Stochastic {PDE}'s and {K}olmogorov equations in infinite
		dimensions ({C}etraro, 1998) 1715\/} (1999), 1--63.
	
	\bibitem{kutyniok2019theoretical}
	{\sc Kutyniok, G., Petersen, P., Raslan, M., and Schneider, R.}
	\newblock A theoretical analysis of deep neural networks and parametric {PDEs}.
	\newblock {\em arXiv:1904.00377\/} (2019), 42 pages.
	
	\bibitem{LabartLelong2013}
	{\sc Labart, C., and Lelong, J.}
	\newblock A parallel algorithm for solving {BSDE}s.
	\newblock {\em Monte Carlo Methods and Applications 19}, 1 (2013), 11--39.
	
	\bibitem{Thomee2003PDEsAndNumerics}
	{\sc Larsson, S., and Thom\'{e}e, V.}
	\newblock {\em Partial differential equations with numerical methods}, vol.~45
	of {\em Texts in Applied Mathematics}.
	\newblock Springer-Verlag, Berlin, 2003.
	
	\bibitem{LeCunBengioHinton15}
	{\sc LeCun, Y., Bengio, Y., and Hinton, G.}
	\newblock Deep learning.
	\newblock {\em Nature 521\/} (2015), 436--444.
	
	\bibitem{LemorGobetWarin2006}
	{\sc Lemor, J.-P., Gobet, E., and Warin, X.}
	\newblock Rate of convergence of an empirical regression method for solving
	generalized backward stochastic differential equations.
	\newblock {\em Bernoulli 12}, 5 (2006), 889--916.
	
	\bibitem{LionnetDosReisSzpruch2015}
	{\sc Lionnet, A., dos Reis, G., and Szpruch, L.}
	\newblock Time discretization of {FBSDE} with polynomial growth drivers and
	reaction-diffusion {PDE}s.
	\newblock {\em The Annals of Applied Probability 25}, 5 (2015), 2563--2625.
	
	\bibitem{LongLuMaDong17}
	{\sc Long, Z., Lu, Y., Ma, X., and Dong, B.}
	\newblock {PDE-N}et: Learning {PDEs} from {D}ata.
	\newblock {\em International Conference on Machine Learning} (2018), 3208--3216.
	
%		\bibitem{LongLuMaDong17}
%	{\sc Long, Z., Lu, Y., Ma, X., and Dong, B.}
%	\newblock {PDE-N}et: Learning {PDEs} from {D}ata.
%	\newblock {\em arXiv:1710.09668\/} (2017), 15 pages.
	
	
	\bibitem{lye2019deep}
	{\sc Lye, K.~O., Mishra, S., and Ray, D.}
	\newblock Deep learning observables in computational fluid dynamics.
	\newblock {\em Journal of Computational Physics}, 109339 (2020), 1--26.
	
%		\bibitem{lye2019deep}
%	{\sc Lye, K.~O., Mishra, S., and Ray, D.}
%	\newblock Deep learning observables in computational fluid dynamics.
%	\newblock {\em arXiv:1903.03040\/} (2019), 57 pages.
	
	\bibitem{MaProtterSanMartin2002}
	{\sc Ma, J., Protter, P., San~Mart\'\i{n}, J., and Torres, S.}
	\newblock Numerical method for backward stochastic differential equations.
	\newblock {\em The Annals of Applied Probability 12}, 1 (2002), 302--316.
	
	\bibitem{MaProtterYong1994}
	{\sc Ma, J., Protter, P., and Yong, J.~M.}
	\newblock Solving forward-backward stochastic differential equations
	explicitly---a four step scheme.
	\newblock {\em Probability Theory and Related Fields 98}, 3 (1994), 339--359.
	
	\bibitem{MaYong1999}
	{\sc Ma, J., and Yong, J.}
	\newblock {\em Forward-backward stochastic differential equations and their
		applications}, vol.~1702 of {\em Lecture Notes in Mathematics}.
	\newblock Springer-Verlag, Berlin, 1999.
	
	\bibitem{magill2018neural}
	{\sc Magill, M., Qureshi, F., and de~Haan, H.}
	\newblock Neural networks trained to solve differential equations learn general
	representations.
	\newblock {\em Advances in Neural Information Processing Systems\/} (2018),
	4075--4085.
	
	\bibitem{Maruyama1955}
	{\sc Maruyama, G.}
	\newblock Continuous {M}arkov processes and stochastic equations.
	\newblock {\em Rendiconti del Circolo Matematico di Palermo. Serie II 4\/}
	(1955), 48--90.
	
	\bibitem{McKean1975}
	{\sc McKean, H.~P.}
	\newblock Application of {B}rownian motion to the equation of
	{K}olmogorov-{P}etrovskii-{P}iskunov.
	\newblock {\em Communications on Pure and Applied Mathematics 28}, 3 (1975),
	323--331.
	
	\bibitem{MilsteinOriginal1974}
	{\sc Milstein, G.~N.}
	\newblock Approximate integration of stochastic differential equations.
	\newblock {\em Theory of Probability \& Its Applications 19}, 3 (1975),
	557--562.
	
	\bibitem{MilsteinTretyakov2006}
	{\sc Milstein, G.~N., and Tretyakov, M.~V.}
	\newblock Numerical algorithms for forward-backward stochastic differential
	equations.
	\newblock {\em SIAM Journal on Scientific Computing 28}, 2 (2006), 561--582.
	
	\bibitem{MilsteinTretyakov2007}
	{\sc Milstein, G.~N., and Tretyakov, M.~V.}
	\newblock Discretization of forward-backward stochastic differential equations
	and related quasi-linear parabolic equations.
	\newblock {\em IMA Journal of Numerical Analysis 27}, 1 (2007), 24--44.
	
	\bibitem{MilsteinTretyakov_SolvingPSPDEsViaAveraging2009}
	{\sc Milstein, G.~N., and Tretyakov, M.~V.}
	\newblock Solving parabolic stochastic partial differential equations via
	averaging over characteristics.
	\newblock {\em Mathematics of Computation 78}, 268 (2009), 2075--2106.
	
	\bibitem{Pham2015}
	{\sc Pham, H.}
	\newblock Feynman-{K}ac representation of fully nonlinear {PDE}s and
	applications.
	\newblock {\em Acta Mathematica Vietnamica 40}, 2 (2015), 255--269.
	
	\bibitem{Raissi18}
	{\sc Raissi, M.}
	\newblock Deep hidden physics models: Deep learning of nonlinear partial
	differential equations.
	\newblock {\em The Journal of Machine Learning Research 19}, 1 (2018), 932--955.
	
	
%		\bibitem{Raissi18}
%	{\sc Raissi, M.}
%	\newblock Deep hidden physics models: Deep learning of nonlinear partial
%	differential equations.
%	\newblock {\em arXiv:1801.06637\/} (2018), 26 pages.
	
	
	\bibitem{RuijterOosterlee2015}
	{\sc Ruijter, M.~J., and Oosterlee, C.~W.}
	\newblock A {F}ourier cosine method for an efficient computation of solutions
	to {BSDE}s.
	\newblock {\em SIAM Journal on Scientific Computing 37}, 2 (2015), A859--A889.
	
	\bibitem{RuijterOosterlee2016}
	{\sc Ruijter, M.~J., and Oosterlee, C.~W.}
	\newblock Numerical {F}ourier method and second-order {T}aylor scheme for
	backward {SDE}s in finance.
	\newblock {\em Applied Numerical Mathematics 103\/} (2016), 1--26.
	
	\bibitem{Ruszczynski2017Dual}
	{\sc Ruszczynski, A., and Yao, J.}
	\newblock A dual method for backward stochastic differential equations with
	application to risk valuation.
	\newblock {\em arXiv:1701.06234\/} (2017), 22 pages.
	
	\bibitem{SirignanoDGM2017}
	{\sc Sirignano, J., and Spiliopoulos, K.}
	\newblock {DGM}: A deep learning algorithm for solving partial differential
	equations.
	\newblock {\em Journal of Computational Physics 375\/} (2018), 1339--1364.
	
	\bibitem{SkorohodBranchingDiffusion1964}
	{\sc Skorokhod, A.~V.}
	\newblock Branching diffusion processes.
	\newblock {\em Theory of Probability \& Its Applications 9}, 3 (1964),
	445--449.
	
	\bibitem{Stroock1982TopicsInSDEs}
	{\sc Stroock, D.~W.}
	\newblock {\em Lectures on topics in stochastic differential equations},
	vol.~68 of {\em Tata Institute of Fundamental Research Lectures on
		Mathematics and Physics}.
	\newblock Tata Institute of Fundamental Research, Bombay; by Springer-Verlag,
	Berlin-New York, 1982.
	\newblock With notes by Satyajit Karmakar.
	
	\bibitem{Thomee1997}
	{\sc Thom\'ee, V.}
	\newblock {\em Galerkin finite element methods for parabolic problems}, vol.~25
	of {\em Springer Series in Computational Mathematics}.
	\newblock Springer-Verlag, Berlin, 1997.
	
	\bibitem{Turkedjiev2015}
	{\sc Turkedjiev, P.}
	\newblock Two algorithms for the discrete time approximation of {M}arkovian
	backward stochastic differential equations under local conditions.
	\newblock {\em Electronic Journal of Probability 20\/} (2015), no. 50, 49.
	
	\bibitem{Watanabe1965Branching}
	{\sc Watanabe, S.}
	\newblock On the branching process for {B}rownian particles with an absorbing
	boundary.
	\newblock {\em Journal of Mathematics of Kyoto University 4\/} (1965),
	385--398.
	
	\bibitem{Zhang2004}
	{\sc Zhang, J.}
	\newblock A numerical scheme for {BSDE}s.
	\newblock {\em The Annals of Applied Probability 14}, 1 (2004), 459--488.
	
\end{thebibliography}

\end{document}